 \documentclass[11pt]{amsart}
\usepackage{amsmath, amsthm, amssymb}
\usepackage{graphicx}
\usepackage{mathbbol}
\usepackage{txfonts}
\usepackage[usenames]{color}

\newtheorem{thm}{Theorem}[section]
\newcommand{\bt}{\begin{thm}}
\newcommand{\et}{\end{thm}}

\newtheorem{conj}[thm]{Conjecture}
\newtheorem{problem}[thm]{Problem}

\newtheorem{cor}[thm]{Corollary}   
\newcommand{\bc}{\begin{cor}}
\newcommand{\ec}{\end{cor}}

\newtheorem{lem}[thm]{Lemma}   
\newcommand{\bl}{\begin{lem}}
\newcommand{\el}{\end{lem}}

\newtheorem{prop}[thm]{Proposition}
\newcommand{\bp}{\begin{prop}}
\newcommand{\ep}{\end{prop}}

\newtheorem{defn}[thm]{Definition}
\newcommand{\bd}{\begin{defn}}    
\newcommand{\ed}{\end{defn}}

\newtheorem{rmrk}[thm]{Remark}   

\newcommand{\br}{\begin{rmrk}}
\newcommand{\er}{\end{rmrk}}

\newcommand{\weaklyto}{\to}

\newtheorem{example}[thm]{Example}

\newcommand{\GHto}{\stackrel { \textrm{GH}}{\longrightarrow} }
\newcommand{\Fto}{\stackrel {\mathcal{F}}{\longrightarrow} }
\newcommand{\VolFto}{\stackrel {V\mathcal{F}}{\longrightarrow} }

\newcommand{\sgn}{\operatorname{sgn}}
\newcommand{\fillvol}{\operatorname{FillVol}}
\newcommand{\mina}{\operatorname{MinA}}
\newcommand{\Sect}{\operatorname{Sect}}
\newcommand{\Ric}{\operatorname{Ric}}
\newcommand{\Scal}{\operatorname{Scalar}}

\newcommand{\be}{\begin{equation}}

 \newcommand{\ee}{\end{equation}}

\newcommand{\N}{\mathbb{N}}

\newcommand{\R}{\mathbb{R}}
\newcommand{\One}{{\bf \rm{1}}}

\newcommand{\E}{\mathbb{E}}

\newcommand{\Z}{\mathbb{Z}}

\newcommand{\diam}{\operatorname{Diam}}
\newcommand{\depth}{\operatorname{Depth}}

\newcommand{\Fm}{{\mathcal F}}

\newcommand{\set}{{\rm{set}}}
\newcommand{\Ricci}{\rm{Ricci}}
\newcommand{\disjointunion}{\sqcup}

\newcommand{\Lip}{\operatorname{Lip}}

\newcommand{\mass}{{\mathbf M}}

\newcommand{\intrectcurr}{{\mathcal I}} 
\newcommand{\intcurr}{{\mathbf I}}      

\newcommand{\area}{\operatorname{Area}}
\newcommand{\vol}{\operatorname{Vol}}

\newcommand{\rstr}{\:\mbox{\rule{0.1ex}{1.2ex}\rule{1.1ex}{0.1ex}}\:}
\begin{document}

\title{Scalar Curvature and Intrinsic Flat Convergence}

\author{Christina Sormani}
\thanks{The author's research was funded by an individual research grant NSF-DMS-1309360 and a PSC-CUNY Research Grant.}

\address{CUNY Graduate Center and Lehman College}
\email{sormanic@gmail.com}




\begin{abstract}
Herein we present open problems and survey examples and theorems concerning sequences of Riemannian manifolds with uniform lower bounds on scalar curvature and their limit spaces.  Examples of Gromov and of Ilmanen which naturally ought to have certain limit spaces do not converge with respect to smooth or Gromov-Hausdorff
convergence.  Thus we focus here on the notion of Intrinsic Flat convergence,
developed jointly with Wenger.  This notion has been applied successfully to study sequences that arise in General Relativity.  Gromov has suggested it should be applied in other settings as well.   We first review intrinsic flat convergence, its properties, and its compactness theorems, before presenting the applications and the open problems.
\end{abstract}

\maketitle

\section{Introduction}
Gromov proved that
sequences of Riemannian manifolds with nonnegative sectional curvature
have subsequences which
converge in the Gromov-Hausdorff sense to Alexandrov spaces with nonnegative Alexandrov curvature \cite{Gromov-metric}.  Burago-Gromov-Perelman
proved that such spaces are rectifiable \cite{BGP}.  
Building upon Gromov's Compactness Theorem,
Cheeger-Colding proved that sequences of Riemannian manifolds 
with nonnegative Ricci curvature 
have subsequences which converge in the metric measure sense to 
metric measure spaces with generalized nonnegative Ricci curvature 
which are also rectifiable \cite{ChCo-PartI} \cite{Gromov-metric}.  

Sequences of manifolds with 
nonnegative scalar curvature need not have subsequences which converge 
in the Gromov-Hausdorff or metric measure sense.  Gromov has suggested that
perhaps under the right conditions
a subsequence will converge in the intrinsic flat sense to
a metric space with generalized nonnegative scalar curvature
\cite{Gromov-Plateau}.  This is an open question: the notion of generalized
scalar curvature has not yet been defined.  

Intrinsic flat convergence was first defined by the author and Wenger in \cite{SorWen2}.  The limits 
obtained under this convergence are countably $\mathcal{H}^m$ rectifiable metric spaces called integral current spaces.   We review the definitions of these notions
within this chapter along with various continuity and
compactness theorems by the author, Perales, Portegies, Matveev, Munn
\cite{Portegies-Sormani}\cite{Perales-Vol}\cite{Perales-Conv}
\cite{Matveev-Portegies}\cite{Munn-F=GH}.   We also review applications of
intrinsic flat convergence to study sequences of manifolds with
nonnegative scalar curvature that arise in General Relativity by the author,
Huang, Jauregui, Lee, LeFloch, and Stavrov 
\cite{Jauregui} \cite{LeeSormani1}\cite{LeFloch-Sormani-1}
\cite{HLS}\cite{Sormani-Stavrov-1}.   We present many examples and state a number of open problems concerning limits of manifolds with
nonnegative scalar curvature.

Recall that a Riemannian manifold, $M^m$, is endowed with a metric tensor,
$g: TM\times TM\to \mathbb{R}$.   One can then define lengths of curves
and distances, 
\be
L(C)=\int_0^1 g(C', C')^{1/2} \, dt,
\,\,\,
d(p,q) = \inf\{ L(C): \, C(0)=p,\, C(1)=q\}.
\ee
If $M$ is compact the distances are achieved as the lengths of
curves called geodesics.   Given any $p\in M$ and any 
vector $V\in T_p(M)$ there is a geodesic,  
\be \label{geod-start}
\gamma(t)=\exp_p(tV)
\textrm{ such that } \gamma(0)=p \textrm{ and } \gamma'(0)=V.
\ee
Taking $e_1...e_m\in TM$ such that  $g(e_i, e_j) =\delta_{i,j}$
one defines Scalar curvature to be the trace of the Ricci curvature
and Ricci to be the trace of the Sectional curvature:
\be \label{Ricci}
\Scal_p= \sum_{i=1}^m \Ric_p(e_i, e_i) \textrm{ where }
\Ric_p(e_i,e_i)= \sum_{j\neq i} \Sect_p(e_i \wedge e_j), 
\ee
\be \label{Sect}
\Sect_p(e_i \wedge e_j)= \lim_{t\to 0}  
6 \left(\frac{tg(e_i,e_j)-d(\exp_p(te_i), \exp_p(te_j))}{t^3}\right).
\ee
Scalar curvature can also be computed using volumes of balls:
\be \label{Scalar-volume}
\Scal_p= \lim_{r\to 0} 6(m+2) 
\left(\frac{\omega_mr^m-\vol(B(p,r))}{\omega_m r^{m+2}}\right).
\ee
In particular,
\be\label{Scalar-volume-0}
\Scal_p >0 0 \iff  \exists r_p>0 \,\,s.t. \,\, \forall r<r_p\,\, \vol(B(p,r)) < \omega_m r^m. 
\ee
This control on volume is too local to apply to prove any global results.  
All properties of manifolds with lower bounds on their
scalar curvature are built using the fact that curvature is defined using tensors 
as in (\ref{Ricci}).   While it may be tempting to
define generalized positive scalar curvature on a limit space using
(\ref{Scalar-volume-0}) it is unlikely to lead to any consequences because the
limit spaces are not smooth and have no tensors.  We need
a stronger definition which implies (\ref{Scalar-volume-0}) and other properties.

Schoen and Yau applied the three
dimensional version of (\ref{Ricci}) to study minimal surfaces in
manifolds with positive scalar curvature.
They proved that a strictly stable closed
minimal surface in a manifold with $\Scal\ge 0$ 
is diffeomorphic to a sphere in
\cite{Schoen-Yau-min-surf}.  In \cite{Schoen-Yau-positive-mass}, they
 applied minimal surface techniques 
to prove the Positive Mass Theorem:
 if $M^3$ is an asymptotically flat Riemannian manifold with
nonnegative scalar curvature then $m_{ADM}(M^3) \ge 0$.
They also proved the following Positive Mass Rigidity Theorem:
\be \label{PMT-Rigidity}
\Scal \ge 0 \textrm{ and }
m_{ADM}(M^3)=0 \implies M^3 \textrm{ is isometric to }\mathbb{E}^3.
\ee
Here $\mathbb{E}^3$ is Euclidean space and the ADM mass is
the limit of the Hawking masses of asymptotically expanding 
spheres
$
m_{ADM}(M) = \lim_{r\to\infty} m_H(\Sigma_r)$ where
\be\label{Hawking-mass}
 m_H(\Sigma)=\sqrt{\frac{\area(\Sigma)}{16\pi}}
\left(1 - \frac{1}{16\pi}\int_\Sigma H^2 \, d\sigma\right).
\ee

Geroch  proved that if $N_t: \mathbb{S}^2 \to M^3$
evolves by inverse mean curvature flow
and $M^3$ has $\Scal \ge 0$ then
then the Hawking mass, $m_H(N_t)$, is nondecreasing. 
Huisken-Ilmanen introduced weak inverse mean curvature flow,
proving it also satisfies Geroch monotonicity and 
$\lim_{t\to\infty} m_H(N_t)= m_{ADM}(M)$.   
They
applied this to prove the Penrose
Inequality: 
\be\label{penrose}
m_{ADM}(M^3)\ge m_H(\partial M^3) = \sqrt{\tfrac{\area(\Sigma)}{16\pi}}
\ee
when $M^3$ is asymptotically flat with a
connected {\em outermost minimizing boundary} (e.g. $\partial M$ is a minimal surface and there are no other closed minimal surfaces in $M$).  
Bray extended their result to have boundaries with more than one connected
component in \cite{Bray-Penrose}.   In addition, there is the Penrose Rigidity Theorem:
\be\label{penrose}
m_{ADM}(M^3)= m_H(\partial M^3) \implies M^3 \textrm{ is isometric to }M_{Sch,m}
\ee
where $M_{Sch,m}$ is the Riemannian Schwarschild space with
mass $m=m_{ADM}(M^3)$.

In addition to Hawking mass, there are other quasilocal masses defined
on manifolds with $\Scal \ge 0$ including the Brown-York mass (which has
nice properties proven by Shi-Tam in \cite{Shi-Tam-2002}) and the Bartnik mass \cite{Bartnik}.  It is not a simple task to define and apply these quasilocal
masses on limit spaces because they all involve the mean curvatures of
surfaces.  Perhaps more promising is Huisken's new isoperimetric quasilocal mass
of a region $\Omega\subset M^3$,
\be\label{mISO}
m_{ISO}(\Omega)=\frac{2}{\area(\partial \Omega)}
\left(\vol(\Omega)-\frac{\area(\partial \Omega)^{3/2}}{6 \sqrt{\pi}}\right),
\ee  
and $m_{ISO}(M)= \limsup_{r\to\infty} m_{ISO}(\Omega_r)$
introduced in \cite{Huisken-isoper-talk}.   Miao has
proven that $m_{ISO}(M)= m_{ADM}(M)$ using volume
estimates of Fan-Shi-Tam in \cite{Fan-Shi-Tam}.
See also work of 
 Jauregui, Lee, Carlotto, Chodosh, and Eichmair 
\cite{Jauregui-Lee:Huisken-isoper} \cite{Carlotto-Chodosh-Eichmair}.   
  
Gromov-Lawson applied the Lichnerowicz formula to prove many things (cf. \cite{Spin-Geometry}) including
the Scalar Torus Rigidity Theorem \cite{Gromov-Lawson-1980} in all dimensions: 
\be\label{Torus-Rigidity}
\Scal \ge 0 \textrm{ and }M^n \textrm{ diffeom to a torus}  \implies 
M^n \textrm{ is isom to a flat torus}.   
\ee
Witten applied 
it to prove the Positive Mass Theorem for Spin manifolds \cite{Witten-PMT}.

Hamilton's Ricci flow leads to a precise control on the scalar curvature 
as well as the areas of minimal surfaces in the evolving manifolds (cf. \cite{Hamilton-survey-95}).  It was applied to prove 
the following rigidity theorems about minimal surfaces.
 Let
\begin{eqnarray}\label{mina}
\mina(M^3)\,\,\, &=& \inf\{ \area(\Sigma^2):
\, \Sigma^2 \textrm{ is a closed min surf in } M^3 \} \textrm{ and }\,\\
\mina_1(M^3) &=& \inf\{ \area(\Sigma^2):
\, \Sigma^2 \textrm{ is a min noncontr sphere in } M^3 \}.
\end{eqnarray}
Note that these invariants are infinite if there are no qualifying minimal surfaces in $M^3$.
Bray, Brendle, and Neves proved the
Cover Splitting
Rigidity Theorem: 
\be \label{Cover-Splitting-Rigidity}
\Scal \ge 2 \textrm{ and }\mina_1(M^3)= 4\pi
\implies \tilde{M}^3 \textrm{ is isom to }\mathbb{S}^2\times \mathbb{R}
\ee
where $\tilde{M}^3$ is the universal cover of $M^3$ \cite{BBN-CAG-10}.
Bray, Brendle, Eichmair and Neves proved the
$\mathbb{RP}^3$ Rigidity
Theorem in \cite{BBEN-CPAM-10}.:  If $M^3$ is diffeomorphic to 
${\mathbb{RP}}^3$ then
\be \label{RP3-Rigidity}
 \Scal \ge 6 \textrm{ and } \mina(M^3) = 2\pi
\implies M^3 \textrm{ is isom to } \mathbb{RP}^3 .
\ee
They have in fact proven theorems which imply these two more simply
stated theorems (\ref{Cover-Splitting-Rigidity})- (\ref{RP3-Rigidity}) as corollaries.

Recall that a rigidity theorem has a statement in the following form: 
\be
M \textrm{ satisfies an hypothesis }  \implies M \textrm{ isometric to } M_0.
\ee
The corresponding almost rigidity theorem (if it exists) would then be:
\be \label{alm-rig-M}
M \textrm{ almost satisfies an hypothesis  }  \implies M \textrm{ is close to } M_0.
\ee
The almost rigidity theorem can also be stated as follows:
\be\label{alm-rig-Mj}
M_j \textrm{ closer and closer to satisfying an hypothesis }
\implies M_j \to M_0.
\ee
Within we describe conjectured almost rigidity theorems for each of the 
rigidity theorems described above.  All of those conjectures remain
open although some have been proven under additional hypothesis.

First one needs to define
closeness for pairs of Riemannian manifolds and convergence of
sequences of Riemannian manifolds.
A pair of compact Riemannian manifolds, $M_1$ and $M_2$, may
be mapped into a common metric space, $Z$, via distance preserving
maps $\varphi_i: M_i\to Z$ which satisfy
\be\label{isom-emb}
d_Z(\varphi_j(x), \varphi_j(y))= d_{M_j}(x,y) \qquad x,y \in M_j.
\ee
Once they lie in a common metric space, $Z$, then one may use the
Hausdorff distance or the flat distance to measure the distance
between the images with respect to the extrinsic space, $Z$.  
We review these extrinsic distances which depend on both
$Z$ and the location of the $M_i$ within $Z$ in Section~\ref{sect-ext}.

However, an intrinsic notion of distance between $M_1$ and $M_2$ can
only depend on intrinsic data about these spaces and not on how they
may be embedded into some extrinsic $Z$.  Thus Gromov defined his
``intrinsic Hausdorff distance'' in \cite{Gromov-metric}, now known as the Gromov-Hausdorff distance,
by taking the infimum over
all distance preserving maps into arbitrary compact metric spaces, $Z$,
of the Hausdorff distance, $d_H^Z$, between the images:
\be
d_{GH}(M_1, M_2) = \inf_{Z, \varphi_i}\left\{ d_H^Z(\varphi_1(M_1), \varphi_2(M_2))
|
\, \varphi_i: M_i\to Z
\right\} 
\ee
Many almost rigidity theorems have been proven for manifolds with
nonnegative Ricci curvature using the Gromov-Hausdorff distance
(cf. \cite{Colding-volume}, \cite{ChCo-almost-rigidity} and \cite{Sor-cosmos}).

For manifolds with nonnegative scalar curvature, one does not
obtain Gromov-Hausdorff closeness in the almost rigidity theorems.
Counterexamples will be described in Section~\ref{sect-ex}.   
Gromov has suggested 
that intrinsic flat convergence may be more well suited towards proving
an Almost Rigidity for the Torus Rigidity Theorem in \cite{Gromov-Dirac}.  Indeed 
some progress has been made by the author, Lee, LeFloch,
Huang, and Stavrov proving special cases of Almost Rigidity for the
Positive Mass Theorem in \cite{LeeSormani1}\cite{LeFloch-Sormani-1}
\cite{HLS}\cite{Sormani-Stavrov-1}.   

The intrinsic flat distance between compact
oriented Riemannian manifolds was defined by the author with
Wenger in \cite{SorWen2} with an infimum over
all distance preserving maps into arbitrary complete metric spaces, $Z$,
of the flat distance, $d_F^Z$, between the images:
\be \label{def-IF-1}
d_{\mathcal{F}}(M_1, M_2) = \inf_{Z, \varphi_i}\left\{ d_F^Z(\varphi_{1\#}[M_1], \varphi_{2\#}[M_2]):
\, \varphi_i: M_i\to Z
\right\}. 
\ee
Intuitively this distance is measuring the filling volume between the two spaces.
One may also consider the intrinsic volume flat distance:
\be\label{def-volF-1}
d_{\vol\mathcal{F}}(M_1, M_2) = d_{\mathcal{F}}(M_1, M_2)
+|\vol(M_1) - \vol(M_2)|
\ee
Full details about the intrinsic flat distance and limits obtained under
intrinsic flat convergence are provided in Section~\ref{sect-int-1}
after a review of Ambrosio-Kirchheim theory in Section~\ref{sect-ext}.


There are a few methods that can be applied to prove almost rigidity theorems.
To apply the {\em explicit control method} one provides enough controls on the $M$
in (\ref{alm-rig-M}) so that one can explicitly construct an embedding of
$M$ and of $M_0$ into a common metric space and explicitly estimate
the distance between them.  This technique was applied to prove GH almost rigidity theorems by Colding in \cite{Colding-volume} and by Cheeger-Colding
in \cite{ChCo-almost-rigidity}.
It was also
applied to prove the $\mathcal{F}$ almost rigidity 
of the Positive Mass Theorem under
additional hypothesis in joint work with Lee \cite{LeeSormani1} and 
in joint work with Stavrov \cite{Sormani-Stavrov-1}.
Lakzian and the author have proven a 
theorem which provides such a construction and estimate if one can
show $M$ and $M_0$ are close on large regions in \cite{Lakzian-Sormani}.
See Section~\ref{sect-cnstr}.

A second technique used to prove almost rigidity theorems is the 
{\em compactness and weak rigidity method}. 
One first provides enough controls on $M_j$
in (\ref{alm-rig-Mj}) so that a subsequence converges to a limit space $M_\infty$.
Then one proves the limit space satisfies the
hypothesis in some weak sense. Finally one proves
the rigidity theorem in that weak setting.   
This technique was applied by the author to prove a GH almost rigidity
theorem in \cite{Sor-cosmos} using Gromov's Compactness Theorem.
which
states that  
\be
Ric_j \ge -(n-1) \textrm{ and }\diam(M_j)\le D \implies M_{j_k} \GHto M_\infty.
\ee
Wenger's Compactness Theorem  
\cite{Wenger-compactness} states that
\be \label{Wenger-compactness-intro}
\diam(M_j)\le D,\,\, \vol(M_j) \le V,\,\, \vol(\partial M_j)\le A \,\,
\implies M_{j_k} \Fto M_\infty.
\ee
Huang, Lee and the author prove Almost Rigidity of the Positive
Mass Theorem for graph manifolds using Wenger's Compactness Theorem combined with an Arzela-Ascoli Theorem and a number of other
theorems concerning intrinsic flat convergence 
in \cite{HLS}.   We will review $\mathcal{F}$
compactness and Arzela-Ascoli theorems in Section~\ref{sect-AA-prop}.

We begin with Section~\ref{sect-ex} surveying examples of sequences
of manifolds with nonnegative scalar curvature.   These examples reveal 
that one cannot simply use intrinsic
flat convergence to handle all the problems that arise when trying to
prove almost rigidity theorems involving nonnegative scalar curvature.
There is a phenomenon called {\em bubbling}.   One may also have
tiny tunnels and construct sequences of manifolds through a process
called {\em sewing} developed by the author with Basilio in \cite{Basilio-Sormani-1},
which lead to limit spaces that do not even satisfy (\ref{Scalar-volume}).
These examples with bubbling and sewing have $\mina(M_j) \to 0$.

We next present the general theory of
Intrinsic Flat convergence and Integral Current Spaces and
survey the key theorems proven in this area.  We begin with
Section~\ref{sect-ext} by reviewing
work of Federer-Flemming and Ambrosio-Kirchheim on integral currents
in Euclidean space and metric spaces. 
In Section~\ref{sect-int-1} we rigorously define Intrinsic Flat Convergence
and Integral Current Spaces and survey known compactness theorems and proposed compactness theorems.   In Section~\ref{sect-cnstr} we present various methods that
may be used to estimate the intrinsic flat distance between two
spaces and describe how these estimates have been used to
prove almost rigidity theorems using the explicit control method.
In Section~\ref{sect-AA-prop} we present theorems about intrinsic flat convergence including theorems about disappearing and
converging points, converging balls, semicontinuity theorems,
Arzela-Ascoli Theorems and Intrinsc Flat Volume Convergence and mention
how these results have been applied to prove almost rigidity theorems using the
compactness and weak rigidity method. 

We close with Section~\ref{sect-conj} which includes statements of conjectures, surveys of partial solutions to the conjectures and recommended related problems.  We discuss the Almost Rigidity of the Positive Mass Theorem, the Bartnik
Conjecture, the Almost Rigidity of the Scalar Torus Theorem, the Almost Rigidity of Rigidity Theorems proven using Ricci Flow, Gromov's Prism Conjecture
and the Regularity of Limit Spaces.   Throughout one hopes to devise a generalized notion of nonnegative scalar curvature on limit spaces. 
Conjectures and problems are interspersed throughout the paper.
If a reader is interested in studying any of these questions, please contact the author.   More details can be provided and the author can coordinate the research of those working on these problems.  

The author must apologize up front that there is no possible way to mention all the
fundementally important papers that have been written concerning manifolds with
scalar curvature bounds.  The results mentioned here have been selected because the author has read them and developed some idea as to how they may applied to
study the intrinsic flat convergence of manifolds with nonnegative scalar curvature.  

\section{Examples with Positive Scalar Curvature} \label{sect-ex}

In this section we survey examples of sequences of three
dimensional Riemannian manifolds, $M_j$, with lower scalar curvature 
bounds and describe their intrinsic
flat limits.  Many of these examples were found by mathematicians
interested in applications to General Relativity.
Manifolds with positive scalar curvature can be viewed as time
symmetric spacelike slices of spacetime satisfying the positive 
energy condition.   Such manifolds are curved by matter and
can have gravity wells and/or black holes with horizons that
are minimal surfaces.   

We do not provide the explicit details or the proofs for these examples
but instead provide references.   We also propose new examples
as open problems that could be written up and published by an interested reader.
We present these examples before presenting intrinsic flat convergence
because they provide some intuitive understanding of intrinsic flat
convergence and what may occur when one has a sequence of
manifolds with nonnegative scalar curvature.   

\subsection{Examples with Wells}\label{sect-ex-wells}

\begin{figure}[htbp]
\begin{center}
\includegraphics[width=8cm]{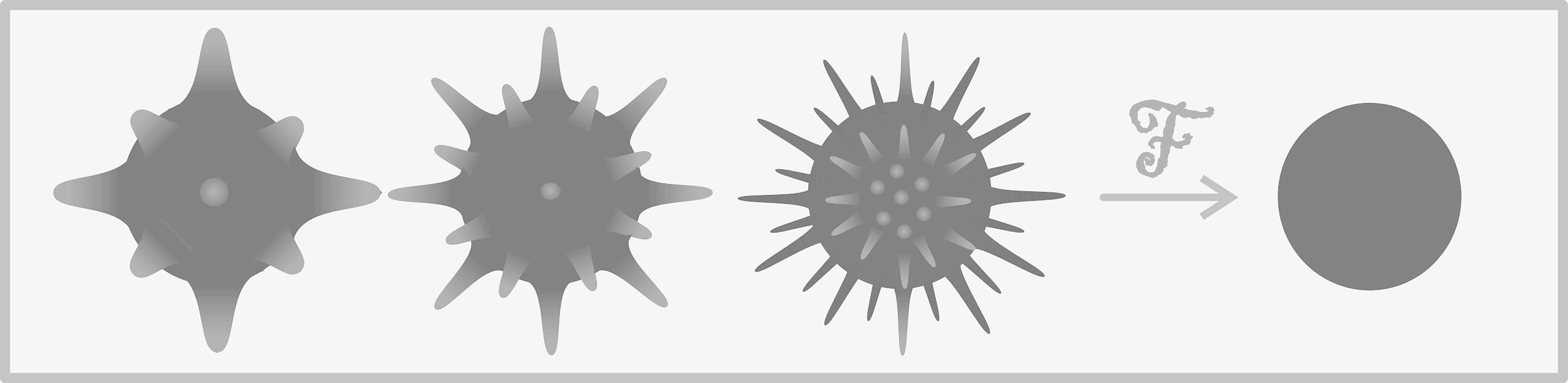}
\caption{The Ilmanen Example.   Image owned by the author.}
\label{fig-Ilmanen}
\end{center}
\end{figure}

Arbitrarily thin arbitrarily deep wells can be constructed 
with positive scalar curvature.  By the Positive Mass 
Theorem, one cannot attach such wells smoothly to Euclidean space.
However they may be glued to spheres of constant positive sectional 
curvature.   In fact, the Ilmanen Example, which initially inspired the
definition of intrinsic flat convergence, consists of a sequence of 
spheres with increasingly many increasingly thin wells as in Figure~\ref{fig-Ilmanen}.   This sequence converges in the intrinsic
flat sense to a standard sphere because a sphere with many thin holes and
a standard sphere can be mapped into a common metric space, and the
the flat distance between them, which is intuitively a filling volume between them,
will be very small.

\vspace{-1cm}

\begin{example} \label{well-sphere}
Lakzian has explicitly
constructed sequences of spheres with one increasingly thin well in 
\cite{Lakzian-Diameter}.  He has also
explicitly constructed the Ilmanen Example of sequences of spheres with
increasingly many increasingly thin wells in the same paper.  He proves
both sequences $\mathcal{F}$ converge to the standard sphere.
He proves the sequence with one increasingly
thin well converges in the GH sense to a sphere with a line segment
attached and the Ilmanen example has no GH limit.
\end{example}

\begin{example} \label{well-mass}
In joint work with Lee, the author
has constructed examples of asymptotically flat
rotationally symmetric manifolds of positive scalar curvature
with arbitrarily thin and arbitrarily deep 
wells and $m_{ADM}(M_j) \to 0$   \cite{LeeSormani1}. 
They have no smooth limits and the
the pointed GH limits of such examples are Euclidean spaces with line
segments of arbitrary length attached.  Lee and the author have also constructed sequences which are not
rotationally symmetric that have $m_{ADM}(M_j) \to 0$ and increasingly
many increasingly thin wells \cite{LeeSormani2}.  Such sequences have no GH converging subsequences because they have increasingly many disjoint balls \cite{Gromov-metric}.   Thus almost rigidity of the positive mass theorem
cannot be proven with GH or smooth convergence, only $\mathcal{F}$
convergence.
\end{example}

\begin{example} \label{well-tori}
It is possible to construct $M_j$ with $Scalar \ge -1/j$ that are diffeomorphic
to tori and contain balls of radius $1/2$ that are isometric to 
balls in rescaled standard spheres.  This will appear in work of the
author with Basilio \cite{Basilio-Sormani-1}.
One may then attach an increasingly thin
well of arbitrary depth to such $M_j$ that have positive scalar curvature.
These examples would $\mathcal{F}$ converge to a standard
flat torus and would GH converge to a standard flat torus with 
a line segment attached.  One may also attach increasingly many
increasingly thin wells of arbitrary depth to the $M_j$ and still $\mathcal{F}$ converge to a standard flat torus but there will be no GH limit of such a sequence.
Thus one must use intrinsic flat convergence to prove almost rigidity
for the Scalar Torus Theorem.
\end{example}

\subsection{Tunnels and Bubbling}

Gromov-Lawson and Schoen-Yau constructed tunnels diffeomoerphic
to ${\mathbb{S}}^2 \times [0,1]$ with
positive scalar curvature which attach 
smoothly on either end to the standard spheres \cite{Gromov-Lawson-tunnels}
\cite{Schoen-Yau-tunnels}.  These tunnels may be arbitrarily thin and long
or thin and short.   At the center of the tunnel, there is a closed minimal surface
diffeomorphic to a sphere.  Sometimes these tunnels are called necks.  

\begin{figure}[h]
\begin{center}
\includegraphics[width=8cm]{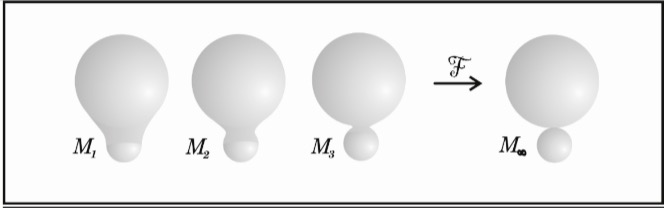}
\caption{Bubbling.    Image owned by the author.}
\label{fig-bubbling}
\end{center}
\end{figure}

\begin{example}\label{tunnel-sphere}
Using these tunnels one may construct sequences of $M_j$
which consist of a pair of standard spheres joined by increasingly thin
tunnels of length $L_j$.  
If $L_j \to 0$, then the GH and $\mathcal{F}$
limit can be shown to be a pair of standard spheres joined at a point
as in Figure~\ref{fig-bubbling}.   This effect is called {\em bubbling}.  
Note that in this example, $\mina(M_j) \to 0$.
\end{example}
\vspace{-.5cm}

\vspace{-.5cm}
\begin{example}\label{ex-not-geod}
If $L_j \to L_\infty>0$, then the GH limit is a pair of standard spheres
joined by a line segment of length $L$ and the $\mathcal{F}$ limit
is just the pair of spheres without the line segment with the restricted 
distance from the GH limit.   Examples similar to these are described
by Wenger and the author in the \cite{SorWen2}.   Notice that the
$\mathcal{F}$ limit is not geodesic.
\end{example}

\begin{example}\label{ex-lots-bubbles}
In fact one may add increasingly many
bubbles with increasingly short and thin tunnels, and the sequence
will have no GH limit and no $\mathcal{F}$ limit.  This does not contradict
Wenger's Compactness Theorem as in (\ref{Wenger-compactness-intro}) 
because the volume is diverging to
infinity even though the diameter is bounded and there is no boundary.
Examples similar to these
with many bubbles of various sizes and tunnels of various lengths
converging to rectifiable limit spaces appear in \cite{SorWen2}.
\end{example}

\begin{example} \label{tunnel-mass}
One may also have bubbling in asymptotically flat sequences of $M_j$
with $m_{ADM}(M_j) \to 0$ obtaining limits which are
Euclidean planes with spheres attached.   This can be done by attaching
bubbles instead of wells to the sequences in Example~\ref{well-mass}.
Such an example demonstrates that any almost rigidity theorem for
the Positive Mass Theorem must somehow avoid bubbling.  In
joint work with Lee, we require that the manifolds have outward minimizing
boundaries just as in the Penrose Inequality \cite{LeeSormani1}.
This is effectively cutting off the bubbles.  One could alternately 
eliminate bubbling
by requiring the sequence to have a uniform lower bound on
the area of the smallest closed minimal surface, $\mina(M_j)\ge A_0>0$. 
\end{example}


\begin{example} \label{tunnel-tori}
One may add a bubble to each $M_j$ 
of Example~\ref{well-tori} with $\Scal_j \ge -1/j$ that are diffeomorphic
to tori and contain balls of constant sectional curvature isometric to
balls in rescaled standard spheres.  Such sequences would $\mathcal{F}$
converge to standard flat tori with a sphere of arbitrary radius
attached at a point.   
This will appear in work of the
author with Basilio \cite{Basilio-Sormani-1}.
One could eliminate such examples
by requiring the sequence have a uniform lower bound on
the area of the smallest closed minimal surface, $\mina(M_j)\ge A_0>0$.   One cannot require 
that there are no closed minimal surfaces here since the
manifolds are diffeomorphic to tori.  
\end{example}

\subsection{Cancellation and Doubling}

The next two examples are described by the author and Wenger
in \cite{SorWen1} \cite{SorWen2}.   Intrinsic flat limit spaces may be
the ${\bf{0}}$ or a rectifiable space with integer weight.
 
\begin{example}\label{cancelling-spheres}
There are sequences of manifolds $M_j^3$ with positive scalar curvature
which have a $\mathcal{F}$ limit which is the
${\bf{0}}$ space,
while converging in the GH sense to a standard three sphere.  This cancelling 
sequence can be constructed with positive scalar
curvature by taking a pair of standard three spheres and connecting them by 
increasingly dense increasingly small tunnels.  These sequences
converge to the ${\bf{0}}$ space because their filling volumes converge to $0$. 
In fact they are the totally geodesic boundaries of four dimensional manifolds whose volume converges to $0$.
\end{example}

\begin{example}
If in the previous example all the tunnels are cut and glued back together with
reversed orientation, then the 
GH limit is still a standard three sphere and
the $\mathcal{F}$ limit is a sphere with weight two everywhere.
\end{example}

\subsection{Sewing Manifolds}

In upcoming joint doctoral work of Basilio with the author, the notion
of sewing Riemannian manifolds is introduced \cite{Basilio-Sormani-1}.  One starts with a 
three dimensional manifold, $M$, that contains a curve, $C:[0,1]\to M$, 
such that a tubular neighborhood around the curve has constant positive
sectional curvature.
One then creates a sequence of manifolds sewn along this curve.  That
is short thin tunnels are attached along the curve pulling the points on the
curve closer together.   The GH and $\mathcal{F}$ limit of such a sequence
is then the original manifold with a {\em pulled thread} along $C$.  That is, all the
points in the image of $C$ have been identified.   One can also sew entire
regions with constant positive sectional curvature to obtain sequences
converging to the original manifold with the entire region identified as a
single point.   If the original manifold has positive scalar curvature, then
so does the sequence.  In addition one may consider sequences of
$M_j$ and sew along curves or in regions of those $M_j$. 
Using this construction, Basilio and the author 
construct the following examples.

\begin{example}\label{sewing-tori}
If one takes the $M_j$ of Example~\ref{well-tori}, one may sew
along curves lying in the balls of radius $1/2$ that have constant
sectional curvature to obtain a sequence of manifolds, $M_j'$ with 
$\Scal_j \ge -1/j$ that are no longer tori but converge to
a limit which is the standard flat torus with a contractible circle 
pulled to a point.   Or a contractible sphere pulled to a point.  Or
a ball of radius $1/2$ pulled to a point.     These examples demonstrate
that limits of manifolds with $\Scal_j\ge -1/j$ may fail to have generalized
nonnegative scalar curvature in the sense that limit in
(\ref{Scalar-volume}) fails to be nonnegative.  
These limits can be biLipschitz to tori and still not be
isometric to a flat tori.  Like all examples created with
this sewing construction, $\mina(M_j)\to 0$.
\end{example}

\vspace{.4cm}
\section{Integral Currents on a Metric Space}\label{sect-ext}

Before we can rigorously 
define intrinsic flat convergence and describe the limit spaces
obtained under intrinsic flat convergence, we need to review 
Ambrosio and Kirchheim's notion of currents and convergence of
currents on a complete metric space \cite{AK}.  Note that like
Federer-Fleming's earlier work on the flat and weak convergence of submanifolds
viewed as currents in Euclidean space, the flat and weak convergence
of Ambrosio-Kirchheim's currents are extrinsic notions of
convergence, depending very much on the way in which the submanifold
or current lies within an extrinsic space.

\subsection{Federer-Fleming currents on Euclidean Space}
In \cite{FF} Federer and Fleming first introduced the
notion of a current on Euclidean space as a generalization of the notion
of an oriented 
submanifold, $\varphi: M \to \E^N$, which views $M$ as a linear functional,
$T=\lbrack M \rbrack$, 
on differential forms:
\be
T(\omega) = [M] \omega =\int_M \varphi^* \omega.
\ee
In particular
\be
T(f \, d\pi_1\wedge \cdots \wedge d\pi_m) =  
\int_M (f\circ \varphi_i)\,
d(\pi_1\circ \varphi_1) \wedge \cdots \wedge d(\pi_m\circ \varphi_m).
\ee
Observe that this is perfectly well defined when $\varphi: M \to \E^N$
is only Lipschitz.   This linear functional captures the notion of boundary,
\be
\partial T(\omega) = \int_{\partial M} \omega = \int_M d\omega=T (d \omega).
\ee
So that
\be
\partial T(f \, d\pi_1\wedge \cdots \wedge d\pi_{m-1})=
T(1 \, df \wedge d\pi_1\wedge \cdots \wedge d\pi_{m-1}).
\ee

Federer and Fleming then studied sequences of submanifolds
by considering the weak limits of their corresponding linear functionals.
They applied this to study the Plateau Problem: searching 
for the submanifold of smallest area with a given boundary.   They proved
sequences of submanifolds approaching the smallest area converge in the
weak sense to a limit which they called an integral current.

\subsection{Ambrosio-Kirchheim Integer Rectifiable Currents}
In \cite{AK}, Ambrosio and Kirchheim defined currents on Euclidean space to integral currents on any complete metric space, $Z$.   
In Federer-Fleming, currents were defined as linear functionals on differential
forms \cite{FF}.   Since there are no differential forms on a metric space, 
Ambrosio and Kirchheim's currents are multilinear functionals which
act on 
DiGeorgi's
$m+1$ tuples \cite{DeGiorgi}.    A {\em tuple} $\left(f,\pi_1, ...,\pi_m\right)$ is in $\mathcal{D}^m(Z)$ iff
$f: Z \to \R$ is a bounded Lipschitz function and
$\pi_i: Z \to \R$ are Lipschitz.  These tuples have no antisymmetry properties.

In \cite{AK} Ambrosio-Kirchheim began their work by defining currents. As we
do not need the notion of a current in this paper.  So we jump directly to 
their notion of an integer rectifiable current applying
their Theorems 9.1 and 9.5 as an explanation rather than using their
definition. 

A linear functional 
$
T: \mathcal{D}^m(Z) \to \mathbb{R}
$
is an $m$ dimensional {\em integer rectifiable current},
denoted $T\in \intrectcurr_m(Z)$ if and only if it can be parametrized
as follows
\be\label{param-representation}
T(f, \pi_1,...,\pi_m) = 
\sum_{i=1}^\infty 
\left(\theta_i \int_{A_i} (f\circ \varphi_i)\,
d(\pi_1\circ \varphi_i) \wedge \cdots \wedge d(\pi_m\circ \varphi_i)\right)
\ee
where $\theta_i \in \Z$ and $\varphi_i:A_i \to \varphi_i(A_i)\subset Z$ are biLipschitz maps defined on 
precompact
Borel measurable sets, $A_i\subset\R^m$, with pairwise disjoint images 
 such that
\be\label{finite-mass}
\sum_{i=1}^\infty\, |\theta_i| \,\mathcal{H}_m(\varphi_i(A_i)) \,\,<\,\, \infty
\quad
\textrm{ where } \mathcal{H}_m
\textrm{ is the Hausdorff measure.}
\ee
A $0$ dimensional integer rectifiable current can be parametrized by
a finite collection of distinct weighted points
\be\label{param-representation-0}
T(f) = 
\sum_{i=1}^N \theta_i f(p_i)
\textrm{ where } \theta_i \in \Z \textrm{ and } p_i \in Z.
\ee

Observe that we then have the following antisymmetry property,
\be
T(f, \pi_1,...,\pi_m)\,= \,\sgn(\sigma) \,T(f, \pi_{\sigma(1)},...,\pi_{\sigma(m)})
\ee
for any permutation $\sigma: \{1,...,m\}\to \{1,...,m\}$.   In addition,
$
T(f, \pi_1,...,\pi_m)\,= 0 
$
if $f$ is the zero function or one of the $\pi_i$ is constant.  So while the tuples do not have the properties of differential forms, the action of the integer rectifiable currents on the tuples has these properties.

Ambrosio-Kirchheim's {\em mass measure} $\|T\| $
of a current $T$, is the smallest Borel measure, $\mu$, such that
\be \label{def-measure-T}
  \Big|T\left(f,\pi\right)\Big|  \,\,  \le \,\,    \int_X |f| d\mu
  \qquad  \forall  \, \left(f,\pi\right) \textrm{ where } \Lip\left(\pi_i\right)\le 1.
\ee  
In Theorem 9.5 of \cite{AK}, the mass measure is explicitly computed.
For the purposes of this paper we need only the following consequence
of their theorem:
\be
m^{-m/2} H_T(A)
\,\,\le\,\, 
||T||(A) 
\,\,\le\,\, 
\frac{2^m}{\omega_m}H_T(A)
\ee
where
\be
H_T(A)=\sum_{i=1}^\infty |\theta_i| \mathcal{H}_m(\varphi_i(A_i)\cap A).
\ee
Furthermore the mass measure of a $0$ dimensional integer rectifiable current 
satisfies 
\be
||T||(A) = \sum_{p_i\in A} |\theta_i|.
\ee

The Ambrosio-Kirchheim {\em mass} of $T$ is defined
\be \label{def-mass-from-current}
M\left(T\right) = || T || \left(Z\right). 
\ee
By the definition of the Ambrosio-Kirchheim mass we have
\be \label{mass-prod-lip}
T(f, \pi_1, ...,\pi_m) \le \sup|f| \,\prod_{i=1}^m \Lip(\pi_i) \,\mass(T).
\ee

The {\em restriction}
of a current $T$ by a $k+1$ tuple 
 $\omega=(g,\tau_1,...\tau_k)\in \mathcal{D}^k(Z)$ with $k<m$ is
 defined by
\be \label{defn-rstr}
(T\rstr\omega)(f,\pi_1,...\pi_m):=T(f\cdot g, \tau_1,...\tau_k, \pi_1,...\pi_m).
\ee
Given a Borel set, $A$, 
\be   
T\rstr A (f,\pi_1,...\pi_m):= T( \One_A \cdot f,\pi_1,...\pi_m)
\ee
where $\One_A$ is the indicator
function of the set.  Observe that $T\rstr \omega$ is an 
integer rectifiable current of dimension $m-k$ and that
\be
\mass(T\rstr \omega) = ||T||(A).
\ee

Given a Lipschitz map, $\varphi:Z\to Z'$, the {\em push
forward} of a current $T$ on $Z$ 
to a current $\varphi_\# T$ on $Z'$is given by
\be \label{def-push-forward}
\varphi_\#T(f,\pi_1,...\pi_m):=T(f\circ \varphi, \pi_1\circ\varphi,...\pi_m\circ\varphi)
\ee
which is clearly still an integer rectifiable current.
Observe that
\be
(\varphi_\#T) \rstr (f, \pi_1,...\pi_k))= \varphi_\#(T \rstr (f\circ \varphi, \pi_1\circ\varphi,...\pi_k\circ\varphi) )
\ee
and
\be \label{rstr-push} 
(\varphi_\#T )\rstr A = (\varphi_\#T) \rstr (\One_A)
=\varphi_\# (T \rstr (\One_A\circ \varphi)) = \varphi_\#(T \rstr \varphi^{-1}(A)).
\ee
In (2.4) of \cite{AK}, Ambrosio-Kirchheim show that
\be  \label{mass-push}
||\varphi_\#T|| \le [\Lip(\varphi)]^m \varphi_\# ||T||,
\ee
so that when $\varphi$ is an isometric
embedding 
\be \label{lem-push-mass}
||\varphi_\#T||=\varphi_\#||T|| \textrm{ and }
\mass(T)=\mass(\varphi_\#T).
\ee

In \cite{AK}[Theorem 4.6] Ambrosio-Kirchheim define the (canonical) set of a current, $T$,
 as the collection of points in $Z$ with positive lower density:
\be \label{def-set-current}
\set\left(T\right)= \{p \in Z: \Theta_{*m}\left( \|T\|, p\right) >0\},
\ee
where the definition of lower density is
\be \label{eqn-lower-density}
\Theta_{*m}\left( \mu, p\right) =\liminf_{r\to 0} \frac{\mu(B_p(r))}{\omega_m r^m}.
\ee
When $T$ is an integer rectifiable current then $\set(T)$ is countably
$\mathcal{H}^m$ rectifiable, which means there exists a collection of
biLipschitz maps, $\varphi_i: A'_i \to \set(T) \subset Z$, defined on Borel sets
$A'_i \in \mathbb{R}^m$ such that
\be
\mathcal{H}_m\left(\set(T) \setminus \bigcup_{i=1}^\infty \varphi_i(A_i) \right)
=0.
\ee
These $\varphi_i$ can be taken from the parametrization of $T$
with $A_i'\subset A_i \subset \bar{A}_i$.

\subsection{Ambrosio-Kirchheim Integral Currents}

The {\em boundary} of $T$ is defined
\be \label{def-boundary}
\partial T(f, \pi_1, ... \pi_{m-1}):= T(1, f, \pi_1,...\pi_{m-1}).
\ee
Note that $\varphi_\#(\partial T)=\partial(\varphi_\#T)$ and 
it can easily be shown that $\partial \partial T=0$. 
The boundary of an integer rectifiable current is not necessarily an integer rectifiable current.
  
An integer rectifiable current  $ T\in\intrectcurr_m(Z)$  is an 
{\em integral current}, denoted $T\in \intcurr_m(Z)$,  if $\partial T$ 
is an integer rectifiable current.  This includes the zero current
\be
0(f,\pi_1,...,\pi_m):=0 \textrm{ with } 
\partial 0(f, \pi_1,...,\pi_{m-1})=0(1,f,\pi_1,...,\pi_{m-1})=0.
\ee
Note that Ambrosio-Kirchheim  
define an integral current
as an integer rectifiable current whose boundary has finite mass
and the more easily applied statement we have here is their Theorem 8.6 in \cite{AK}.

Given an
oriented Riemannian manifold with boundary, $M^m$, 
such that $\vol_m(M)<\infty$ and $\vol_{m-1}(\partial M) < \infty$, and
given a Lipschitz
map $\varphi: M \to Z$, we can define an integral current
$\varphi_\#[M]\in \intcurr_m(Z)$ as follows
\be \label{subman}
\varphi_\#[M] (f, \pi_1,..., \pi_m)=\int_M (f\circ \varphi)\, 
d(\pi_1\circ \varphi) \wedge \cdots \wedge d(\pi_m\circ \varphi).
\ee 
Note that $\partial \varphi_\#[M] = \varphi_\#[\partial M] $ where
$\partial M$ is the boundary of $M$ and 
\be
\mass(\varphi_\#[M])=\vol_m(\varphi(M)).
\ee
If $\vol_m(M)<\infty$ and $\vol_{m-1}(\partial M) =\infty$, then
$[M]$ is only integer rectifiable and not integral.


Whenever, $T$ is an integral current, $\partial \partial T=0$, and
\be
\partial: \intcurr_m(Z) \to \intcurr_{m-1}(Z).
\ee
In addition, if 
$\varphi: Z_1 \to Z_2$ is Lipschitz, then by (\ref{def-push-forward})
\be
\varphi_{\#}: \intcurr_m(Z_1) \to \intcurr_{m}(Z_2).
\ee
The restriction of an integral current defined in (\ref{defn-rstr}) need not be 
an integral current.  However,
the Ambrosio-Kirchheim Slicing Theorem implies that
\be\label{rstr-balls}
T\rstr B(p,r) \textrm{ is an integral current for almost every }r>0
\ee  
where $B(p,r)=\{x:\, d(x,p)<r\}$.  


\subsection{Convergence of Currents in a Metric Space}
 
In Definition 3.6 of \cite{AK}, Ambrosio and Kirchheim state that a
sequence of integral currents $T_j \in \intcurr_m\left(Z\right)$ 
lying in a complete metric space, $Z$,
is said to converge weakly to a current $T$, denoted $T_j \weaklyto T$,
 iff the pointwise limits satisfy
\be
\lim_{j\to \infty}  T_j\left(f, \pi_1,...\pi_m\right) = T\left(f, \pi_1,...\pi_m\right) 
\ee
for all bounded Lipschitz $f: Z \to \R$ and Lipschitz $\pi_i: Z \to \R$.
Ambrosio-Kirchheim next observe that if $T_j$ converges weakly to $T$, then 
the boundaries converge
\be\label{bndry-conv}
\partial T_j \weaklyto \partial T,   
\ee
and the mass is lower semicontinuous
\be\label{semicont}
\liminf_{j\to\infty} \mass(T_j) \ge \mass(T).
\ee
Thus the weak limit of a sequence of integer rectifiable
currents with a uniform upper bound on mass is an
integer rectifiable current:
\be
T_j \in \intrectcurr_m(Z), \,\,\,\mass(T_j) \le V_0
\textrm{ and }
T_j \to T \implies T \in \intrectcurr_m(Z).
\ee
Similarly for integral currents we have
\be
T_j \in \intcurr_m(Z), \,\,\,\mass(T_j) \le V_0,\,\,\,\mass(\partial T_j) \le A_0
\textrm{ and }
T_j \to T \implies T \in \intcurr_m(Z).
\ee
For any open set, $A\subset Z$, if $T_j \to T$ then
\be
\liminf_{j\to\infty} ||T_j||(A) \ge ||T||(A).  
\ee
However $T_j \rstr A$ need not converge weakly to $T \rstr A$
(cf. Example 2.21 of \cite{Sormani-AA}).

Ambrosio-Kirchheim prove the following compactness theorem:

\begin{thm}\label{AK-compact}\cite{AK}
Given any complete metric space 
$Z$, a compact set $K \subset Z$ and $A_0, V_0>0$.
Given
any sequence of integral currents  $T_j \in \intcurr_m \left(Z\right)$ satisfying
\be\label{AK-compactness}
\mass(T_j) \le V_0 \textrm{, } \mass(\partial T_j) \le A_0
\textrm{ and }
\set\left(T_j\right) \subset K,
\ee there exists a subsequence, $T_{j_i}$, 
which converges weakly to $T\in \intcurr_m(Z)$.
\end{thm}

It is possible that the limit obtained in this theorem
is the $0$ integral current.  Observe that whenever the sequence of
currents is {\em collapsing}, 
\be
\mass(T_j) \to 0,
\ee
then by (\ref{mass-prod-lip}) we have
\be
|T_j(f, \pi_1, ...,\pi_m)| \le \sup|f| \,\prod_{i=1}^m \Lip(\pi_i) \,\mass(T_j) \to 0
\ee
and so $T_j$ converges weakly to $0$.

It is also possible for $T_j$ to converge weakly to $0$ without collapsing.
This can occur due to {\em cancellation}, when the $T_j$ fold over on themselves
as in Example~\ref{ex-taco}.  We include this 
example in detail because it inspires the notion of flat convergence and will
be refered to repeatedly in this paper.

\begin{example} \label{ex-taco}
Let $T_j = \varphi_{j\#}[M]  \subset \intcurr_2(\E^3)$ where
\be \label{folding}
\varphi_j(s,t)=\left(s, \,t/j, \,|t|b_j/j\right) \textrm{ where } b_j=\sqrt{j^2-1}
\ee
on  $M=\{(s,t):\, s\in [-1,1],\, t\in [-1,1]\}$.
Since 
$
\mass(\varphi_{j\#}[M])=\vol(M) 
$
does not converge to $0$, this sequence is not collapsing.
Observe that 
$
T_j = A_j + \partial B_j
$ where
\be \label{Bj}
B_j=[\{(x,y,z): \, |x|\le 1,\, |y|\le 1/j, \, z\in [|y|b_j, b_j/j] \, \}]
\in \intcurr_3(\E^3) \textrm{ and}
\ee
\be \label{Aj}
A_j\,\,\,= A_{j}^-+A_{j}^+ +A_{j}^0 \in \intcurr_2(\E^3)\textrm{ where}
\ee
\begin{eqnarray}
A_j^- &=& -[\{\, (-1,y,z):\, \, y\in [-1/j, 1/j], \, z\in [|y| b_j, b_j/j \,]\, \}\rbrack  \\
A_j^+ &=& [\{\, (+1,y,z):\, \, y\in [-1/j, 1/j], \, z\in [|y| b_j, b_j/j \,]\, \}\rbrack   \\
A_j^0 &=& [\{\, (x, y, b_j/j\,):\,  x\in [-1,1],\, y\in [-1/j, 1/j]\, \} \rbrack.
\end{eqnarray}
Since
\be\label{AjBjmass}
\mass(B_j)\le (4/j)\textrm{ and }
\mass(A_j) \le  (2/j) + (2/ j) + 4/j,
\ee
we have
$
B_j \weaklyto 0 \textrm{ and } A_j \weaklyto 0.
$
By (\ref{bndry-conv}) we have $\partial B_j \weaklyto \partial 0=0$, and thus
\be
T_j=A_j +\partial B_j \weaklyto 0.
\ee
\end{example}

Sometimes part of a sequence disappears under weak convergence
and part remains.  This happens in the following example:

\begin{example} \label{one-spline-to-disk}
Let $T_j = \varphi_{j\#}[D^2]  \subset \intcurr_2(\E^3)$ with
$
D^2=\{(x,y):\,\, x^2+y^2 \le 1\}
$
and
\be
\varphi_{j\#}(x,y)=(x,y, f_j(\sqrt{x^2+y^2})),
\ee
where $f_j:[0,1]\to [0,1]$ is a smooth cutoff function such that $f_j(r)=1$
near $r=0$ and $f_j(r)=0$ for $r\ge 1/j$.  Then 
$\partial T_j = \varphi_{j\#}[S^1]$ is constant and so the sequence 
does not disappear.   In fact $T_j$ converges weakly to
$T_\infty=\varphi_{\infty\#}[D^2]$ where
\be
\varphi_{j\#}(x,y)=(x,y, 0)
\ee
since $T_j-T_\infty= \partial B_j$ where
\be \label{Bj2}
B_j= [\{(x,y,z): \, x^2+y^2 \le 1,\,\, 0\le z \le f_j(\sqrt{x^2+y^2})\}].
\ee
Since 
\be \label{Bjmass}
\mass(B_j) \le \pi (1/j)^2 \to 0
\ee
we have $B_j \to 0$ and thus $T_j - T_\infty \to 0$ and $T_j \to T_\infty$.
\end{example}

\subsection{The Flat distance vs the Hausdorff distance}\label{subsect-H}

In \cite{Wenger-flat}, Wenger defines the flat distance between 
two integral currents, $T_1, T_2\in \intcurr_m(Z)$, lying in a common complete
metric space, $Z$, to be
\be \label{Flat-in-Z}
d_F^Z\left(T_1,T_2\right)
= \inf\left\{ \mass(A) +\mass(B):\,\, A+\partial B= T_1-T_2 \right\}
\ee
where the infimum is taken over all $A \in \intcurr_m(Z)$ and $B \in \intcurr_{m+1}(Z)$
such that $A+\partial B= T_1-T_2$.  This is the same definition 
given by
Federer and Fleming in \cite{FF} building on work of Whitney \cite{Whitney}
for the flat distance in Euclidean
space, where it is a norm, $|T_1-T_2|_\flat$.   The lack of scaling in
(\ref{Flat-in-Z}) is a result of setting the flat distance to be a norm on Euclidean
space.   A scalable version of the flat distance might be defined for
$Z$ with a finite diameter $\diam(Z)=D$ as follows (cf. \cite{LeFloch-Sormani-1}):
\be \label{Flat-in-Z}
d_{DF}^Z\left(T_1,T_2\right)
= \inf\left\{ D \mass(A) +\mass(B):\,\, A+\partial B= T_1-T_2 \right\}
\ee

Observe that if two oriented hypersurfaces share a boundary, then the flat distance between them is
bounded above by the volume between them.  In Example~\ref{one-spline-to-disk}, we have
\be
d_F^{\E^3}(T_j, T_\infty) \to 0
\ee
by taking $B_j$ as in (\ref{Bj2}) 
so that $T_j-T_\infty=\partial B_j $.  Then
$d_F^Z(T_j , 0) \le \mass(B_j)$
which converges to $0$ as $j \to \infty$ by (\ref{Bjmass}).
In Example~\ref{ex-taco}, we
produce a sequence of integral currents $T_j$ in Euclidean space such that
\be
d_F^{\E^3}(T_j, 0) \to 0
\ee
by taking $B_j$ as in (\ref{Bj}) and $A_j$ as in (\ref{Aj})
so that $T_j-0=\partial B_j + A_j$.  Then
$d_F^Z(T_j , 0) \le \mass(A_j) +\mass(B_j)$
which converges to $0$ as $j \to \infty$ by (\ref{AjBjmass}).

In \cite{Wenger-flat}, Wenger proves that
when
\be
\mass(T_j)\le V_0
\textrm{ and }
\mass(\partial T) \le A_0
\ee
then weak and flat convergence are equivalent:
\be \label{F=weak}
T_j \weaklyto T
\textrm{ if and only if }
d_F^Z(T_j, T) \to 0. 
\ee

One should contrast the flat distance between submanifolds
viewed as integral currents with the Hausdorff distance between
submanifolds viewed as subsets, $X_i=\varphi_i(M_i)$.  Note
that the Hausdorff distance is defined to be
\be \label{Hausdorff}
d_H^Z(X_1, X_2) =
\inf\{ r>0: \,\, X_1\subset T_r(X_2),\,
 X_2 \subset T_r(X_1)\}
 \ee
 where $T_r(X)= \{z:\, \exists x \in X \, s.t.\, d(x,z)<r\}\subset Z$.
There is no notion of
 a disappearing Hausdorff limit.  The Hausdorff limit of a collapsing 
 sequence of sets like $[0,1/j]\times [0,1] \subset \mathbb{E}^2$
 is easily seen to be $\{0\}\times [0,1] \subset \mathbb{E}^2$, which 
 is simply a lower dimensional set.   The Hausdorff limit of
 the sequence of cancelling submanifolds, $\varphi_j(M)$, in 
 Example~\ref{ex-taco} 
 is easily seen to be the set $[-1,1]\times\{0\}\times [0,1]$.   
No points in a Hausdorff limit can disappear.   In Example~\ref{one-spline-to-disk}, 
the Hausdorff limit of $\varphi_j(D^2)$ is 
a disk with a line segment attached:
$
(D^2\times \{0\}) \cup (\{0,0\}\times [0,1]) \subset {\mathbb{E}}^3.
$

\begin{figure}[htbp]
\begin{center}
\includegraphics[width=4in]{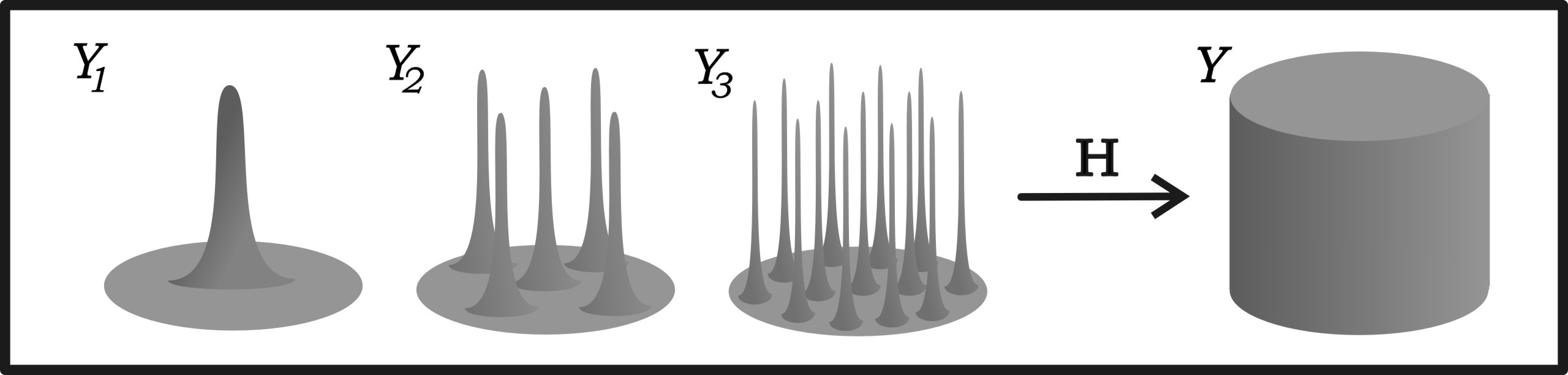}
\caption{A troublesome Hausdorff limit.  Image owned by the author.}
\label{fig-Plateau}
\end{center}
\end{figure}

One reason Federer and Fleming introduced integral currents and
flat convergence was to solve the Plateau problem of finding
a minimal surface with a given boundary, $\Gamma$.  Suppose
for example that
\be
\Gamma=\{(\cos(t), \sin(t),0): \, t\in \mathbb{S}^1\}\subset\E^3. 
\ee
One must find the surface 
$\varphi(D^2)$ such that $\partial(\varphi(D^2))=\Gamma$
of smallest area.   One may try to find this minimal surface by taking
a sequence of such surfaces, $\varphi_j(D^2)$, with area decreasing to
the infimum of these areas, and look for a limit.  In Example~\ref{one-spline-to-disk}, we have such a sequence with a thinner and thinner spine
so that the Hausdorff limit is a disk with a line segment attached,
not a minimal surface.   Even worse, one may
have a sequence of $\varphi_j(D^2)=Y_j$ with increasingly many increasingly dense spines 
as in Figure~\ref{fig-Plateau}
so that the Hausdorff
limit is 
\be
Y= \{ (x,y,z): \, x^2+y^2 \le 1, \,\, z\in [0,1]\}.
\ee
  This Hausdorff limit has no notion of boundary and
is no longer two dimensional.  There is no smooth or even $C_0$ limit of
such $\varphi_j$.

On the other hand
the flat limit of a sequence,
$\varphi_j(D^2)$, with area decreasing to
the infimum of these areas does exist and
is the standard disk
\be
\varphi_\infty: \{(x,y): \, x^2+y^2 \le 1\} \to \E^3
\textrm{ with } \varphi_\infty(x,y)=(x,y,0).
\ee
This was proven in the case with one spline in Example~\ref{one-spline-to-disk}.
This can be seen in Figure~\ref{fig-Plateau} because the volume between the
$\varphi_j(M_j)$ and $\varphi_\infty(M_\infty)$ converges to $0$.  

Even on a compact metric space, $Z$, flat convergence is well suited to
the Plateau problem, where one is given $\Gamma \in \intcurr_{m-1}(Z)$
and asked to find an integral current $T\in \intcurr_m(Z)$, such that
$\partial T = \Gamma$ and
\be
\mass(T) =  M_0=\inf\{\mass(T): \, \partial T=\Gamma\}.
\ee
One takes $T_j \in \intcurr_m(Z)$ such that $\partial T_j = \Gamma$,
and $\mass(T_j) \to M_0$.  By Ambrosio Kirchheim's Compactness
Theorem [Theorem~\ref{AK-compact} above] and (\ref{F=weak}), 
a subsequence converges in the weak and flat sense to some
$T_\infty \in \intcurr_m(Z)$.
Since $\partial T_\infty = \Gamma$
by (\ref{bndry-conv}) and $\mass(T_\infty)=M_0$
by (\ref{semicont}), we have a desired solution to the Plateau problem.

\vspace{.4cm}
\section{Integral Current Spaces and Intrinsic Flat Convergence}\label{sect-int-1}

In this section we provide the rigorous definition for the
intrinsic flat convergence of a sequence of oriented Riemannian
manifolds or, more generally, a sequence of integral current spaces
\cite{SorWen2}.   It is crucial to remember that the manifolds in the
sequence are not submanifolds of any common Euclidean space.  This
is an intrinsic notion about the intrinsic geometry of the Riemannian manifolds.

As described in the introduction, the intrinsic flat distance is
defined much like the Gromov-Hausdorff distance, by taking an
infimum over all distance preserving maps, $\varphi_j: M^m_j \to Z$
into any common complete metric
space, $Z$:
\be
d_{\mathcal{F}}(M_1, M_2) = \inf_{Z, \varphi_j}\left\{ d_F^Z\left(\varphi_{1\#}[M_1], \varphi_{2\#}[M_2]\right)\right\}
\ee
where this is now rigorously defined using (\ref{subman}) and (\ref{Flat-in-Z}).
In fact we can define the intrinsic flat distance between $M_1$ and
an abstract $\bf{0}$ space as well:
\be
d_{\mathcal{F}}(M_1, {\bf{0}}) = \inf_{Z, \varphi_j}\left\{ d_F^Z\left(\varphi_{1\#}[M_1], 0\right)\right\}.
\ee
Keep in mind that by the Kuratowski Embedding Theorem
any pair of separable metric spaces can be isometrically embedded
into a Banach space, $Z$, so these infima are always finite (cf. \cite{SorWen2}).
In \cite{SorWen2}, the author and Wenger prove that for $M_j$ compact, we have
$d_{\mathcal{F}}(M_1, M_2)=0$ if and only if there is an orientation
preserving isometry between them.

It is essential to remember that the $\varphi_j$ are distance preserving maps
or isometric embeddings in the sense of Gromov as in (\ref{isom-emb}).  
They are not Riemannian
isometric embeddings which only preserve lengths of curves.  For example,
the Riemannian isometry from the standard circle, ${\mathbb{S}}^1$, to the boundary of the closed Euclidean
disk, $D^2$, is not a distance preserving map.  The Riemannian isometry from the standard circle, $\mathbb{S}^1$, to the boundary of the heimsphere,
${\mathbb{S}}^2_+$ is a distance preserving map.   
In Example~\ref{ex-taco} we have a single flat square,
$M=[-1,1]\times [-1,1]$, with a sequence of $\varphi_j: M \to \mathbb{E}^3$
which preserve lengths of curves, and yet the flat limit of the images
is $0$ due to cancellation.   If the intrinsic flat distance were defined using
such maps, then  $d_{\mathcal{F}}(M_1, {\bf{0}})=0$, and similarly the
intrinsic flat distance between any pair of oriented manifolds would be $0$.

In this section we introduce a larger class of spaces, integral current spaces,
which are metric spaces
with an additional structure.   These spaces include oriented Riemannian
manifolds with boundary and their intrinsic flat limits.   We then define the intrinsic flat distance between this larger class of spaces and review fundamental theorems
about intrinsic flat convergence.   

\subsection{Integral Current Spaces}\label{subsect-int-cur-space}

Unlike the Gromov-Hausdorff distance, the intrinsic flat distance 
cannot be defined between an arbitrary pair of metric spaces, $M_j=(X_j,d_j)$.
One needs an additional structure which guarantees that the isometric
embeddings of the $M_j$ into $Z$ may be viewed as integral currents.
Thus the author and Wenger introduced the following notion
in \cite{SorWen2}:

\begin{defn} \label{defn-int-curr-space}
An $m$ dimensional metric space 
$M=\left(X,d,T\right)$ is called an integral current space if
it has a integral current structure $T \in \intcurr_m\left(\bar{X}\right)$
where $\bar{X}$ is the metric completion of $X$
and $\set(T)=X$.   Also included in the $m$ dimensional integral
current spaces is the ${\bf 0}$ space, denoted ${\bf{0}}=(\emptyset,0,0)$.
We say two such spaces are equal, $M_1=M_2$, if there is a
current preserving isometry, $F: M_1 \to M_2$:
\be \label{current-preserving-isometry}
d_2(F(x), F(t))=d_1(x,y) \textrm{ and } F_\# T_1= T_2.
\ee
The mass of the integral current space is, $\mass(M)=\mass(T)$.
The diameter is $\diam(M)=\sup\{d(x,y):\, x,y\in X\}$ if $M\neq 0$
and $\diam({\bf{0}})=0$. 
\end{defn}

Any $m$ dimensional integral current space $(X,d,T)$ is countably 
$\mathcal{H}^m$ rectifiable in the sense that there
exists biLipschitz charts $\varphi_i: A_i \to X$ 
where $A_i\subset \R^m$ are Borel
and
\be \label{lim-rect}
\mathcal{H}^m\left(X\setminus \bigcup_{i=1}^\infty \varphi_i(A_i) \right)=0.
\ee
In fact these charts can be viewed as oriented with
weights $\theta_i\in \N$ 
as in (\ref{param-representation}). Any $0$ dimensional integral current space,
$(X,d,T)$ is a finite collection of points, $X=\{p_1,...,p_N\}$,
with a metric $d$ and with 
weights $\theta_i \in \mathbb{Z}$ so that
$T(f)= \sum_{i=1}^N \theta_i f(p_i)$.  
  
A compact oriented Riemannian manifold with boundary, $(M^m, g)$,
is an integral current space, where $X=M^m$, where $d=d_g$ is the standard
metric on $M$,
\be\label{dg}
d_g(p,q) = \inf \{L_g(C):\, C(0)=p, \, C(1)=q\}
\ee
where
\be
 L_g(C)=\int_0^1 g(C'(t), C'(t))^{1/2} \, dt,
\ee
and where $T=[M]$ is integration over $M$,
\be \label{manifold}
T(f, \pi_1,...,\pi_m)= \int_M f\, d\pi_1 \wedge \cdots \wedge d\pi_m.
\ee
In this setting $\mass(M)=\vol(M)$.   Note that if $(M_1, g_1)$
and $(M_2, g_2)$ are diffeomorphic then they have the same integral
current structure up to a sign.  In fact they need only be biLipschitz equivalent.  They do not have the
same mass unless there is a volume preserving diffeomorphism
between them.  They are not viewed as the same integral current space unless there is an orientation preserving isometry between them.

If $M$ is a precompact oriented Riemannian manifold with boundary,
$(M^m, g)$, then we can define an integral current space $(X,d,T)$,
by taking the metric completion 
$\bar{X}=\bar{M}$, defining $d_g$ as the continuous
extension to $\bar{X}$ of (\ref{dg}), and defining $T \in \intcurr_m(\bar{X})$
exactly as in (\ref{manifold}).   Then $X=\set(T) \subset \bar{X}$.
This set is called the {\em settled completion} of $M$ and is denoted, $M'$.
 In particular
 \be \label{settled}
 M' = \{ x\in \bar{M}: \,\, \liminf_{r\to 0} \vol_m(B(x,r)\cap M)/r^m >0 \}.
 \ee 
 So if $M$ is a manifold with a singular point removed, that point is
 always included in the metric completion $\bar{M}$ but it is not included
 in the settled completion if it is a cusp singularity.

The boundary of an integral current space, $\left(X,d,T\right)$,
is the integral current space:
\be
\partial \left(X,d_X,T\right) := \left(\set\left(\partial T\right), d_{\bar{X}}, \partial T\right)
\ee
where the distance on the boundary is $d_{\bar{X}}$ which is restricted
from the distance on the metric completion $\bar{X}$.
If $\partial T=0$ then one says $\left(X,d,T\right)$ is an integral current without boundary.  The ${\bf 0}$ space has no boundary.

Note that the boundary of $(M, d_g, [M])$ when $M$ is an oriented Riemannian manifold with boundary is
$(\partial M, d_g, \partial[M])$ endowed with the restricted distance, $d_g$,
defined on $M$ as in (\ref{dg}).  It is only a geodesic space if $\partial M$
is totally geodesic in $M$.   For example, when $M=D^2\subset \E^2$, then
$\partial M=(\mathbb{S}^1, d_{\E^2}, [\mathbb{S}^1])$ is not a geodesic space because
there are no curves whose length is equal to the distance between the 
points, 
\be
d_{\E^2}(p,q)=|p-q|< d_{\mathbb{S}^1}(p,q)=\cos^{-1}( 1-|p-q|^2/2).   
\ee
The boundary of the upper hemisphere, 
$M= \left(\mathbb{S}^2_+, d_{\mathbb{S}^2_+}, [\mathbb{S}^2_+]\right)$,
is a geodesic integral current space, 
$\partial M=\left(\mathbb{S}^1, d_{\mathbb{S}^1}, [\mathbb{S}^1]\right)$.




In \cite{Sormani-AA} the author proves that
a ball in an integral current space, $M=\left(X,d,T\right)$,
with the current restricted from the current structure of $M$ is an integral current space itself, 
\be\label{Spr}
S(p,r):=\left(\,\set(T\rstr B(p,r)),d,T\rstr B\left(p,r\right)\right)
\ee
for almost every $r > 0$.   Furthermore,
\be\label{ball-in-ball}
B(p,r) \subset \set(S(p,r))\subset \bar{B}(p,r)\subset X.
\ee
Note that
the outside of the ball, $(M\setminus B(p,r), d, T-S(p,r))$, 
and the sphere, 
\be
\partial S(p,r):=\left(\,\set(\partial(T\rstr B(p,r))),d,\partial(T\rstr B\left(p,r\right))\right),
\ee
are integral current spaces 
for the same values of $r>0$.  If $\partial M=0$ then
\be
\set(\partial(T\rstr B(p,r))) \subset \{x: \, d(x,p)=r\}.
\ee

In \cite{Portegies-Sormani}, Portegies and Sormani investigate
the notion of the filling volume
\be \label{def-fillvol}
\fillvol(\partial M)=\inf\{ \mass(N): \,\, \partial N=\partial M\}
\ee
where the infimum is over all integral current spaces, $N$, such that there is a current preserving isometry from 
$\partial N$ to $\partial M$.   This notion of filling volume does not quite agree
with Gromov's notion of Filling Volume in \cite{Gromov-filling} because in our
notion
there is a larger collection of candidates, $N$, for filling the manifold, because we require a current preserving isometry on the boundary, and because we use the Ambrosio-Kirchheim mass.   With our notion, we immediately have
\be
\mass(M) \ge \fillvol(\partial M).
\ee

Applying the filling volume to balls, we have for almost every $r>0$ that
\be
||T||(B(p,r))=\mass(T\rstr B(p,r))=\mass(S(p,r)) \ge \fillvol(\partial S(p,r)).
\ee
In particular, if a point $x\in \bar{X}$ satisfies
\be
\liminf_{r\to 0} \fillvol(\partial S(p,r))/r^m >0
\ee
then $x\in \set(T)=X$.   This idea was first applied jointly with
Wenger in \cite{SorWen1}
before the notion of integral current space was precisely defined in \cite{SorWen2}.
Further exploration of filling volume and a new notion called 
the sliced filling volume appeared in \cite{Portegies-Sormani} with Portegies.

\vspace{.4cm}
\subsection{Intrinsic Flat Convergence} \label{subsect-int-flat}

The intrinsic flat distance between integral current spaces
was first defined by the author and Wenger in \cite{SorWen2}:

\begin{defn} \label{def-flat1} 
 For $M_1=\left(X_1,d_1,T_1\right)$ and $M_2=\left(X_2,d_2,T_2\right)\in\mathcal M^m$ let the 
intrinsic flat distance  be defined:
 \begin{equation}\label{equation:def-abstract-flat-distance}
  d_{\Fm}\left(M_1,M_2\right):=
 \inf d_F^Z
\left(\varphi_{1\#} T_1, \varphi_{2\#} T_2 \right),
 \end{equation}
where the infimum is taken over all complete metric spaces 
$\left(Z,d\right)$ and distance preserving maps 
$\varphi_j : \left(\bar{X}_j,d_j\right)\to \left(Z,d\right)$. 
\end{defn}

When $M_j$ are precompact integral current spaces
 we prove the infimum in this definition
is obtained \cite{SorWen2}[Thm 3.23] and consequently
$d_\mathcal{F}$ is a distance \cite{SorWen2}[Thm3.27] on the class of precompact integral current spaces up to
current preserving isometries as in (\ref{current-preserving-isometry}).
In particular, it is a distance on the class of oriented compact manifolds with boundary of a given dimension.

We say 
\be
M_j \Fto M_\infty \textrm{ iff } d_{\mathcal{F}}(M_j, M_\infty) \to 0.
\ee
By the definition $M_j \Fto M_\infty$ if and only if 
there exists distance preserving maps to complete metric spaces, 
$\varphi_j: M_j \to Z_j$ and $\varphi'_j: M_\infty \to Z_j$, and integral
currents, $B_j \in \intcurr_{m+1}(Z_j)$ and $A_j \in \intcurr_{m}(Z_j)$,
such that 
\be
\varphi_{j\#} T_j -\varphi'_{j\#} T_\infty= \partial B_j + A_j
\ee
and
\be
d_{\mathcal{F}}(M_j, M_\infty)\le d_F^{Z_j}
\left(\varphi_{j\#} T_j, \varphi'_{j\#} T_\infty \right)\le 
\mass(B_j)+\mass(A_j) \to 0.
\ee
We could then replace $Z_j$ with $Z'_j$ that are closures of countably
$\mathcal{H}^{m+1}$ rectifiable spaces by taking 
\be \label{Zj-here}
Z'_j=Cl(\set(B_j) \cup \set(A_j)) \subset Z_j.
\ee
So in fact the $Z$ in the infimum of the Definition~\ref{def-flat1} may
be chosen in this class.

Note that if $M_j \Fto M_\infty$ then using the same distance preserving
maps we have
\be \label{bndry-1}
\varphi_{j\#} \partial T_j -\varphi'_{j\#} \partial T_\infty= \partial A_j
\ee
and
\be \label{bndry-2}
d_{\mathcal{F}}(\partial M_j, \partial M_\infty)\le d_F^{Z_j}
\left(\varphi_{j\#} \partial T_j, \varphi'_{j\#} \partial T_\infty \right)\le 
\mass(A_j) \to 0.
\ee
So $\partial M_j \Fto \partial M_\infty$.

The following theorem in \cite{SorWen2} is an immediate consequence
of Gromov's Embedding Theorem and Ambrosio-Kirchheim's Compactness Theorems:

\begin{thm} \label{GH-to-flat}
Given a sequence of precompact $m$ dimensional integral current spaces $M_j=\left(X_j, d_j, T_j\right)$ such that 
\be
\left(\bar{X}_{j}, d_{j}\right) \GHto \left(Y,d_Y\right),\,\,\,
\mass(M_j)\le V_0 \,\,\,\textrm{ and } \,\,\, \mass(\partial M_j)\le A_0
\ee
then a subsequence converges 
 in the 
intrinsic flat sense 
\be
\left(X_{j_i}, d_{j_i}, T_{j_i}\right) \Fto \left(X,d_X,T\right)
\ee
where either $\left(X,d_X,T\right)$ is the ${\bf 0}$ integral current space
or $\left(X,d_X,T\right)$ is an $m$ dimensional integral current space
such that  $X \subset Y$ with the restricted metric $d_X=d_Y$.
\end{thm}

Immediately one notes that if $Y$ has Hausdorff dimension less than $m$,
then $(X,d,T)=\bf{0}$.   In Section~\ref{subsect-compactness} we survey theorems
in which it is proven under additional hypothesis that the intrinsic
flat and GH limits agree.   There are many examples
with nonnegative scalar curvature where these limits do not agree 
presented in \cite{SorWen2}. In fact one may not even have a GH converging subsequence
for a sequence with an intrinsic flat limit (see Section~\ref{sect-ex-wells}).

Gromov's Embedding Theorem which is applied to prove Theorem~\ref{GH-to-flat}, 
states that if $\left({X}_{j}, d_{j}\right) \GHto \left(X_\infty,d_\infty\right)$ then
there is a compact metric space $Z$ and a collection of isometric
embeddings $\varphi_j: X_j \to Z$ such that 
\be \label{Gromov-Z}
d_H^Z(\varphi_j(X_j), \varphi_\infty(X_\infty)) \to 0.
\ee 
Note that without his embedding theorem one needs different $Z_j$ for
each term in the sequence and then one would not be able to apply the
Ambrosio-Kirchheim Compactness Theorem 
(cf. Theorem~\ref{AK-compact}) to complete the proof of 
Theorem~\ref{GH-to-flat}.   

In \cite{SorWen2}[Thms 4.2-4.3], the author and Wenger prove similar
embedding theorems for sequences which converge in the intrinsic flat sense:
if 
\be
M_{j}=\left(X_j, d_j, T_j\right) \Fto M_0=\left(X_\infty,d_\infty,T_\infty\right),
\ee 
then
there is a common separable
complete metric space, $Z$, and distance preserving maps  $\varphi_j: X_j \to Z$ such that
\be\label{converge}
d_F^Z(\varphi_{j\#}T_j,\varphi_{\infty\#} T_\infty)\to 0.
\ee
In the case where $M_0=\bf{0}$ then we have (\ref{converge})
as well with $\varphi_{\infty\#} T_\infty=0$ and can find $z\in Z$
and $x_j \in X_j$ such that $\varphi_j(x_j)=z$.   In fact this $Z$ can be chosen
to be the closure of a countably $\mathcal{H}^{m+1}$ rectifiable metric space and is
glued together from the $Z'_j$ in (\ref{Zj-here}).

These embedding theorems
 do not
 require uniform bounds on the masses or volumes of the $M_j$ and $\partial M_j$.
 Combining them with Ambrosio-Kirchheim's lower semicontinuity
 of mass (\ref{semicont}) we see that
 \be \label{semi-mass}
M_j \Fto M_\infty \implies \liminf_{i\to\infty} \mass(M_i) \ge \mass(M_\infty).
 \ee
In \cite{Sormani-AA} the author proves lower semicontinuity of the
diameter as well:
\be \label{semi-diam}
M_j \Fto M_\infty \implies \liminf_{i\to\infty} \diam(M_i) \ge \diam(M_\infty),
 \ee
 
In \cite{Portegies-Sormani}, the author and Portegies prove
that
\be \label{fillvol-cont}
\partial M_j \Fto \partial M_\infty \implies \fillvol(\partial M_j) \to \fillvol(\partial M_\infty).
\ee
This idea was first observed in joint work of the author with Wenger \cite{SorWen1}.
Since
\be
\mass(M) \ge \fillvol(\partial M)
\ee
one can use filling volumes to provide a lower bound on the 
mass of the limit
\be
\mass(M_\infty) \ge \fillvol(\partial M_\infty) =\lim_{j\to \infty} \fillvol(\partial M_j).
\ee
 Portegies and the author also introduce the notion of a {\em sliced filling
 volume} in \cite{Portegies-Sormani} and prove that it is continuous
 with respect to intrinsic flat convergence and also provides a
 lower bound for mass.
 
\subsection{Compactness Theorems for Intrinsic Flat Convergence}
\label{subsect-compactness}

The first compactness theorem for intrinsic flat convergence is stated by the
author with Wenger in
\cite{SorWen2}.  It is
a combination of Gromov's Compactness and Embedding Theorems 
with Ambrosio-Kirchheim
Compactness, to say that if $M_j=(X_j, d_j, T_j)$ satisfy the hypothesis
of Gromov's Compactness Theorem and of Ambrosio-Kirchheim's Compactness
theorem, then a subsequence converges in the GH sense and the $\mathcal{F}$ sense
where the $\mathcal{F}$ limit is a subset of the GH limit 
(cf. Theorem~\ref{GH-to-flat}).  There are
a number of theorems which apply Gromov's Compactness
theorem combined with this theorem and then prove
the GH and $\mathcal{F}$ limits agree under additional hypothesis
including noncollapsing, $\mass(M_j) \ge V_0>0$.  
We call these $\mathcal{F}=$GH compactness theorems.   

The author and Wenger prove in
\cite{SorWen2} a $\mathcal{F}=$GH compactness theorem
for sequences of manifolds without boundary that either
have uniform linear contractibility functions or 
are noncollapsing with $\Ric \ge 0$.  Perales has extended this to allow 
boundaries with
various conditions on the boundary in \cite{Perales-Conv}.   Matveev-Portegies
have extended the result without boundary to uniform negative lower
bounds on Ricci curvature in \cite{Matveev-Portegies}.   The author,
Huang and Lee have proven a $\mathcal{F}=$GH compactness
theorem for sequences of integral current spaces, $(X, d_j, T)$, with 
varying bounded distance functions $d_j$ in the Appendix to \cite{HLS}.
Li and Perales have proven a $\mathcal{F}=$GH compactness theorem for 
noncollapsing integral current
spaces $(X_j, d_j, T_j)$ with nonnegative Alexandrov
curvature (including manifolds with nonnegative sectional curvature)
in \cite{Li-Perales}.    
It is unknown whether integral current
spaces satisfying various generalized notions of Ricci curvature
have $\mathcal{F}=$GH compactness theorems.

In the setting with $\Scal \ge 0$, we do not in general have GH limits and
so we need compactness theorems with weaker hypothesis that do not
imply GH convergence of subsequences. 
Wenger's Compactness Theorem was proven in \cite{Wenger-compactness}
and stated in the following form in \cite{SorWen2}:

\begin{thm}{\em \bf Wenger Compactness}\\
If $M_j$ are integral current spaces of dimension $m$ satisfying the following
\be\label{Wenger-comp}
\mass(M_j) \le V_0 \qquad \mass(\partial M_j) \le A_0 \qquad \diam(M_j) \le D_0
\ee
then there exists a subsequence $M_{j_k} \Fto M_\infty$ where
$M_\infty$ is an integral current spaces of dimension $m$
possibly $\bf{0}$.
\end{thm}

Perales has applied this theorem in \cite{Perales-Vol}
to prove two $\mathcal{F}$ compactness
theorems.  One assumes the given sequence of oriented manifolds satisfies
\be\label{PeralesV1}
\Ric_j \ge 0 \qquad \vol(\partial M_j) \le A_0 \qquad H_{\partial M_j} \ge H_0
\qquad \diam(M_j) \le D_0
\ee
and the other assumes the given sequence satisfies
\be\label{PeralesV2}
\Ric_j \ge 0 \qquad \vol(\partial M_j) \le A_0 \qquad H_{\partial M_j} \ge H_0>0
\qquad \diam(\partial M_j) \le D_0.
\ee
Here $H$ is the mean curvature with respect to the outward pointing normal.
Note that in (\ref{PeralesV2}) the only condition on the interior of the
manifold is $\Ric\ge 0$.

Observe that in both of these theorems, we could renormalize the
manifolds to have $\vol(\partial M_j) = A_0$.   When $H_0\ge 0$, these sequences
have Hawking mass as in (\ref{Hawking-mass}) uniformly bounded above:
\be
 m_H(\partial M_j)
\le m_0=\sqrt{\frac{A_0}{16\pi}}
\left(1 - \frac{1}{16\pi}A_0 H_0^2\right).
\ee
This leads naturally to the following conjecture 
which could be a step towards proving almost rigidity of
the Positive Mass Theorem or Bartnik's Conjecture \cite{Bartnik}:

\begin{conj}\label{Hawking-Compactness}{\em \bf Hawking Mass Compactness}\\
Given a sequence of 
three dimensional oriented manifolds $M_j^3$ satisfying 
\be
 \vol(M_j) \le V_0\qquad 
\vol(\partial M_j) = A_0 
\qquad \diam(M_j) \le D_0.
\ee
\be\label{scalar-hawking-bounds}
\qquad\,\, \Scal_j \ge 0 \qquad \qquad H_{\partial M_j} \ge 0 \qquad
\qquad m_H(\partial M_j)\le m_0 \qquad
\ee
and either no closed interior minimal surfaces or $\mina(M_j) \ge A_1>0$,
then a subsequence converges in the intrinsic flat sense
\be
M_{j_k} \Fto M_\infty \textrm{ and } \mass(M_{j_k}) \to \mass(M_\infty)
\ee
and $M_\infty$ satisfies (\ref{scalar-hawking-bounds})
in some generalized sense (cf. Section~\ref{sect-generalized}).
One might replace Hawking mass with another quasilocal mass in this conjecture.
\end{conj}

LeFloch and the author have proven this 
Hawking Mass Compactness Conjecture in the rotationally
symmetric setting assuming that there are no closed interior minimal
surfaces in \cite{LeFloch-Sormani-1}.  This is shown by proving $H^1_{loc}$
convergence of a subsequence of the manifolds with a well chosen
gauge and then proving the $H^1_{loc}$ limit is a $\mathcal{F}$ limit
using Theorem~\ref{thm-subdiffeo}.
In general it is unknown whether $H^1_{loc}$ convergence implies
$\mathcal{F}$ convergence, but here there is also monotonicity of the
Hawking mass to help.   Since the limit space is a rotationally symmetric
manifold with a metric tensor $g \in H^1_{loc}$, it is possible to define
generalized notions of nonnegative scalar curvature using (\ref{Ricci})
in a weak sense and also to define Hawking mass and show 
(\ref{scalar-hawking-bounds}) hold on the limit spaces as well.   

Gromov has conjectured vaguely that intrinsic flat convergence
may preserve some notion of $\Scal \ge 0$ in \cite{Gromov-Plateau} and
\cite{Gromov-Dirac}.   Considering the
examples and the above conjecture, we propose the following Scalar
Compactness Theorem which requires a uniform lower bound on
the area of a closed minimal surface: 
 
\begin{conj}\label{Scalar-Compactness} {\em \bf Scalar Compactness}\\
Given a sequence of 
 oriented manifolds $M_j^3$ with $\partial M^3_j = {\bf{0}}$
satisfying
\be
 \vol(M_j) \le V_0\quad \diam(M_j) \le D_0\quad
\Scal_j \ge 0 \quad \mina(M_j) \ge A_1>0
\ee
then a subsequence converges in the intrinsic flat sense
\be
M_{j_k} \Fto M_\infty \textrm{ and } \mass(M_{j_k}) \to \mass(M_\infty)
\ee
and $M_\infty$ has $\Scal_\infty\ge 0$
in some generalized sense (cf. Section~\ref{sect-generalized}). 
\end{conj} 

A proof of this Scalar Compactness Theorem 
in the rotationally symmetric case might imitate the proof of
the Hawking Compactness Theorem
of the author with LeFloch, however the work in \cite{LeFloch-Sormani-1} very strongly uses that there is a boundary to choose a gauge.  Nevertheless
a very similar proof should work and would make a nice problem
for a doctoral student.   Quite a different technique would be needed
to handle other settings.  In the graph case 
as in \cite{HLS} there is no $H^1_{loc}$ convergence.

See Section~\ref{sect-generalized} for more about how
generalized $\Scal_\infty\ge 0$ might be defined.

\section{Theorems which imply Intrinsic Flat Convergence} \label{sect-cnstr}
 
In this section we present theorems which have been applied to prove 
sequences of spaces converge in the intrinsic flat sense.  In the previous section
we have already presented compactness theorems which imply intrinsic
flat convergence of subsequences.   
Here we present theorems where geometric constraints
and relationships between a pair of spaces are used to bound
the intrinsic flat distance between them.

Before we begin, note that in Section 5 of \cite{SorWen2}, the author and Wenger
proved that if the Gromov Lipschitz distance between two Riemannian
manifolds is small, then the intrinsic flat distance is small.  In particular
if the manifolds are close in the $C_0$ sense then they are close in the
intrinsic flat sense.  This theorem may now be viewed as a special
case of Theorem~\ref{thm-subdiffeo} included below.   More general
statements about pairs of integral current spaces with such bounds
also appear in \cite{SorWen2}.
 
\subsection{Using Riemannian Embeddings to estimate $D_{\mathcal{F}}$} 
 
 Recall that in the definition of intrinsic flat convergence, one
 must find distance preserving maps of the pair of manifolds
 into a common complete metric space,
 $Z$, before estimating the flat distance between the images.   
 If, however, one only has Riemannian isometric embeddings of
 the manifolds into a common Riemannian manifold, then one 
 may apply the following theorem proven by Lee and the author
 in \cite{LeeSormani1} to estimate the intrinsic flat distance between the
 spaces.
  
  \begin{thm}\label{embed-const}
If $\varphi_i: M^m_i \to N^{m+1}$ are Riemannian isometric embeddings with 
embedding constants $C_{M_i}$ where
\be \label{eqn-embed-const-1}
C_M:= \sup_{p,q\in M} | d_M(p, q) - d_N(\varphi(p),\varphi(q)) |,
\ee
and if
they are disjoint and lie in the boundary of a region $W \subset N$
then
\begin{eqnarray}
d_{\mathcal{F}}(M_1, M_2) &\le& 
S_{M_1}\left(\vol_m(M_1)+ \vol_{m-1}(\partial M_1) \right) \\
&&+S_{M_2}\left(\vol_m(M_2)+ \vol_{m-1}(\partial M_2) \right)\\
&&+ \vol_{m+1}(W) + \vol_{m}(V)
\end{eqnarray}
where $V= \partial W \setminus \left( \varphi_1(M_1) \cup \varphi_2(M_2)\right)$
where $S_{M_i}= \sqrt{C_{M_i}(\diam(M_i)+C_{M_i})}$.
\end{thm}

This theorem is proven in \cite{LeeSormani1} by explicitly constructing
a geodesic metric space,
\be
Z= W_0 \cup W_1 \cup W_2 
\subset N \times [0, S_M],
\ee
where $W_0=\left\{(x,0): \,x\in N\right\}$ and
\be
W_i=\left\{(x,s): \, x\in \varphi_i(M_i),\, s\in [0, S_{M_i}]\right\}
\ee
and proving $\psi_i(x)=(\varphi_i(x), S_{M_i})$ are distance preserving maps into $Z$.  Then taking $B=[W_0]+[W_1]-[W_2]$ and $A=V_0 + V_1 - V_2$ where
$V_0=[ \{(x,0): \,\, x\in V\}]$ and
\be
V_i=[\left\{(x,s): \, x\in \varphi_i(\partial M_i),\, s\in [0, S_{M_i}]\right\}].
\ee
with the appropriate orientations, we get $\psi_{1\#}[M_1]-\psi_{2\#}[M_2]=A +\partial B$. The estimate then follows because $d_{\mathcal{F}}(M_1, M_2)\le \mass(A) + \mass(B)$.
 
\subsection{Smooth Convergence Away from Singular Sets} \label{sect-smooth}

If one has a pair of Riemannian manifolds containing subregions that
are close in the $C_0$ sense, then one can estimate the intrinsic
flat distance between these manifolds using estimates on the
volumes of the regions where they are different and additional information
as proven by Lakzian and the author in \cite{Lakzian-Sormani}: 

\begin{thm} \label{thm-subdiffeo}
Suppose $M_1=(M,g_1)$ and $M_2=(M,g_2)$ are oriented
precompact Riemannian manifolds
with diffeomorphic subregions $U_i \subset M_i$ and
diffeomorphisms $\psi_i: U \to U_i$ such that the following hold
\be \label{thm-subdiffeo-1}
(1+\epsilon)^{-2} \psi_2^*g_2(V,V)<\psi_1^*g_1(V,V)
< (1+\epsilon)^2 \psi_2^*g_2(V,V) \qquad \forall \, V \in TU,
\ee
\be \label{DU}
D_{U_i}= \sup\{\diam_{M_i}(W): \, W\textrm{ is a connected component of } U_i\},
\ee
\be \label{lambda}
\lambda=\sup_{x,y \in U}
|d_{M_1}(\psi_1(x),\psi_1(y))-d_{M_2}(\psi_2(x),\psi_2(y))|.
\ee
Then the Gromov-Hausdorff distance between the metric
completions is bounded,
\be \label{thm-subdiffeo-6}
d_{GH}(\bar{M}_1, \bar{M}_2 ) \le a + 2\bar{h} +
\max\left\{ d^{M_1}_H(U_1, M_1), d^{M_2}_H(U_2, M_2)\right\}
\ee
and the intrinsic flat distance between the
settled completions is bounded,
\begin{eqnarray}
d_{\mathcal{F}}(M'_1, M'_2) &\le&
\vol_m(M_1\setminus U_1)+\vol_m(M_2\setminus U_2) \label{sda1}\\
&&+
\left(2\bar{h} + a\right) \Big(
\vol_{m-1}(\partial U_{1})+\vol_{m-1}(\partial U_{2})\Big)\label{sda2}\\
&&+\left(2\bar{h} + a\right) \Big(
\vol_m(U_{1})+\vol_m(U_2))\Big) \label{sdb}.
\end {eqnarray}
where
$
a=a(\epsilon, D_{U_1}, D_{U_2}) 
$
converges to $0$ as $\epsilon \to 0$ for fixed values of $D_{U_i}$,
and where
$
h = h(\epsilon, \lambda, D_{U_1}, D_{U_2})
$
converges to $0$ as both $\epsilon \to 0$ and $\lambda\to 0$
for fixed $D_{U_i}$.  Explicit formulas for $a(\epsilon, D_{U_1}, D_{U_2})$
and $h(\epsilon, \lambda, D_{U_1}, D_{U_2})$ are given in \cite{Lakzian-Sormani}.
\end{thm}

\begin{figure}[htbp]
\begin{center}
\includegraphics[width=4in]{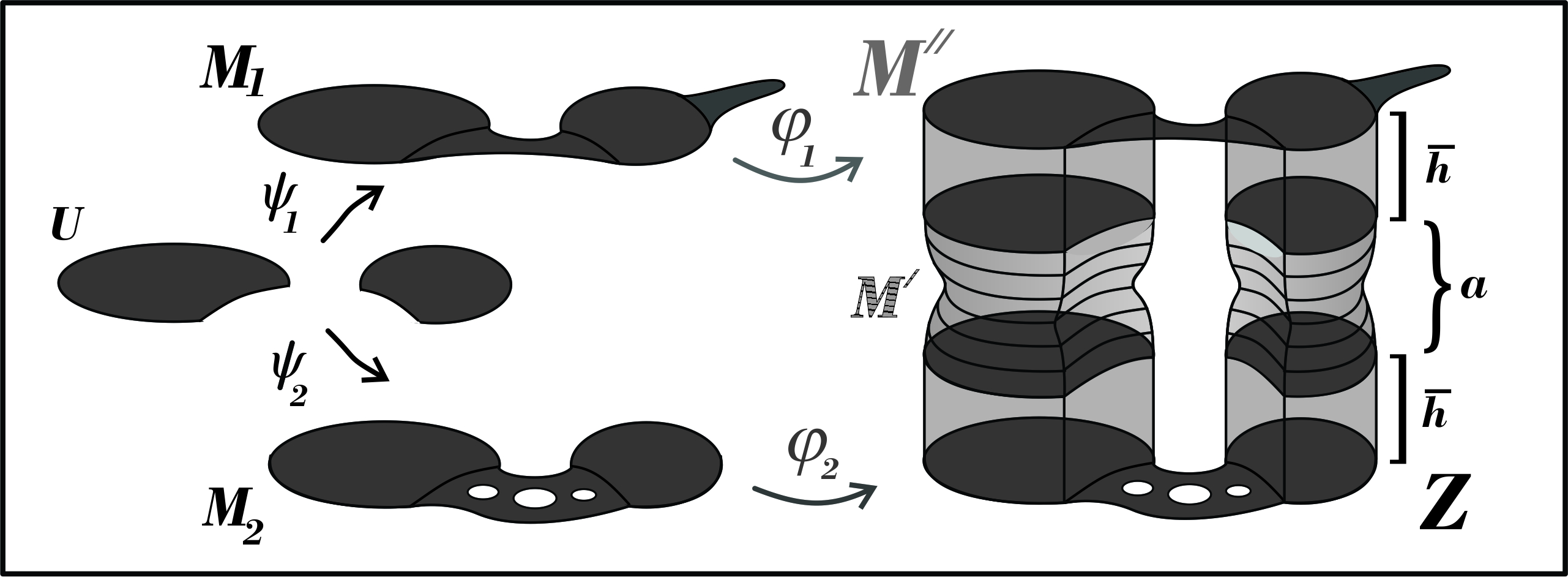}
\caption{Estimating the intrinsic flat distance.    Image owned by the author
and Sajjad Lakzian.}
\label{fig-LS}
\end{center}
\end{figure}

Theorem~\ref{thm-subdiffeo} is proven in \cite{Lakzian-Sormani}
by explicitly constructing a 
common metric space 
\be
Z = M_1\disjointunion_{U_1}(U\times [0,2\bar{h}+a] ) \disjointunion_{U_2}M_2
\ee
where $U\times [0,2\bar{h}+a]$ is a product manifold with a precisely given
metric $g'$, and
$M_1$ is glued to it along $U_1\subset M_1$ which is
isometric to $U \times\{0\}$
and  $M_2$ is glued to it along $U_2\subset M_2$ which is
isometric to $U \times\{2\bar{h}+a\}$ as in Figure~\ref{fig-LS}.
The metric $g'$ is chosen so that $\varphi_i: M_i \to Z$
are distance preserving.  In particular $g'=dt^2+ g_1$ on
$U\times [0,\bar{h}]$ where $\bar{h}$ is chosen as in the theorem
statement to guarantee that there are no
short paths between points in $M_1$ that run through $U\times [0,2\bar{h}+a]$.
Similarly, $g'=dt^2+ g_2$ on $U\times [a+\bar{h}, a+2\bar{h}]$.   On the
middle, $U\times [\bar{h}, a+\bar{h}]$, there is a hemispheric warping
between the metrics $g_1$ and $g_2$.   Once an explicit $Z$ has
been found, then an explicit $A\in \intcurr_m(Z)$ is found where
$\mass(A)$ is the sum of the terms in (\ref{sda1})-(\ref{sda2}) and
an explicit $B\in \intcurr_{m+1}(Z)$ is found where $\mass(B)$ is
the term in (\ref{sdb}).   The details are easy to follow in
\cite{Lakzian-Sormani}.

Theorem~\ref{thm-subdiffeo}
has been applied in work of the author and Lakzian
to prove the rotationally symmetric Hawking Mass Compactness Theorem
\cite{LeFloch-Sormani-1}.  It has been applied by Lakzian to study
Ricci flow through neck pinch singularities in \cite{Lakzian-Cont-Ricci}.
It has been applied by Lakzian in \cite{Lakzian-Diameter} to prove 
$\mathcal{F}$ convergence of sequences of manifolds which converge 
smoothly away from singular sets.   In particular, Lakzian proves that if 
$M_j=(M,g_j)$ be a sequence of compact oriented Riemannian manifolds
with a set $S\subset M$ such that $\mathcal{H}^{n-1}(S)=0$,
and a connected precompact exhaustion,
$W_k$, of $M\setminus S$ satisfying 
\be \label{defn-precompact-exhaustion}
\bar{W}_k \subset W_{k+1} \textrm{ with } 
\bigcup_{k=1}^\infty W_k=M\setminus S
\ee
with $g_j$ converge smoothly to $g_\infty$ on each $W_k$,
\be\label{m-diam}
\diam_{M_j}(W_k) \le D_0 \qquad \forall i\ge j, 
\ee
\be \label{m-area}
\vol_{g_j}(\partial W_k) \le A_0,
\ee
and
\be \label{m-edge-volume}
\vol_{g_j}(M\setminus W_k) \le V_k \textrm{ where } \lim_{k\to\infty}V_k=0.
\ee
Then
\be
\lim_{j\to \infty} d_{\mathcal{F}}(M_j', N')=0.
\ee
where  $N'$ is the settled completion
of $N=(M\setminus S, g_\infty)$ as in (\ref{settled}).
If one also assumes $\Ric_j\ge H$, then 
\be
\lim_{j\to \infty} d_{GH}(M_j, \bar{N})=0.
\ee
Other theorems about smooth convergence away from singular
sets are proven in \cite{Lakzian-Diameter} as well.   KAHLER

\subsection{Pairs of Integral Current Spaces}

The following theorem concerns pairs of integral current spaces
which share the same space, $X$, and the same current structure, $T$,
but have different distance functions.  This happens for example when
one has a pair of oriented Riemannian manifolds with an orientation
preserving diffeomorphism between them but can also be applied to
pairs of integral current spaces with a biLipschitz current preserving map
between them.
The Gromov-Hausdorff part of this theorem 
was proven by Gromov in \cite{Gromov-metric} and the intrinsic
flat part was proven by the author with Huang and Lee
in the Appendix to \cite{HLS}. 

\begin{thm}\label{app-thm}
Fix a precompact $m$-dimensional integral current space $(X, d_0, T)$
with $\partial T=0$ and fix
$\lambda>0$.   Suppose that
$d_1$ is another distance on $X$ such that
\be\label{eps}
\epsilon= \sup\left\{|d_1(p,q)-d_0(p,q)|:\,\, p,q\in X\right\}. 
\ee 
Then we have the following:
\be
d_{GH}\left((X, d_1), (X, d_0)\right) \le 2\epsilon
\ee
\be 
d_{\mathcal{F}}\left((X, d_1, T), (X, d_0, T)\right) \le
2^{(m+1)/2} \lambda^{m+1} (2\epsilon) \mass(T).
\ee
\end{thm}

Note that the hypothesis of this theorem arises when
an integral current space, $(X,d,T)$, has an almost distance
preserving map into another
metric space, $F: X \to Y$, and one defines a new distance
$d_1$ on $X$ as $d_1(p,q)=d_Y(F(p), F(q))$.   It was applied by the author
with Huang and Lee in \cite{HLS}
as one of the final steps towards proving the almost rigidity for the 
Positive Mass Theorem in the graph setting.

\begin{problem}
It seems natural that one could apply Theorem~\ref{app-thm} to
prove a more general theorem about pairs of integral current
spaces, $(X_i, d_i, T_i)$, which contain a pair of regions
$U_i \subset X_i$ with a current preserving isometry
between the regions.  If one simply uses the restricted 
distances then there is no need for an embedding constant
or a $\lambda$ as in the prior two sections.  However, one
needs to extend the Theorem~\ref{app-thm} to consider
the setting with boundary.
\end{problem}

\section{Theorems about Intrinsic Flat Limits}\label{sect-AA-prop}
 
In this section we assume we have a sequence, 
$M_j \GHto M_\infty$ or $M_j \Fto M_\infty$
and present theorems about limits of points in these spaces,
limits of functions on these spaces, continuity and semicontinuity of
various quantities on these spaces.   Recall that we
have already mentioned a number of such results in 
Subsection~\ref{subsect-int-flat} and we will be refering to
those results here as well.   We close this section with a discussion
of the setting where one has intrinsic flat convergence with
volume continuity.
 
\subsection{Limits of Points and Points with no Limits}  

When studying sequences of converging Riemannian manifolds, $M_j$,
one often wishes to understand what happens to points, $x_j \in M_j$: {\em where
do they converge and when do they disappear?}
One is also interested in understanding limit points, $x_\infty \in M_\infty$,
by considering $x_j\in M_j$ that converge to these points.  This can
easily be understood when the $M_j$ are a sequence of converging
submanifolds lying in a given space, and it is clear that under flat
convergence of submanifolds some sequences of points will disappear
in the limit.  Defining converging sequences of points when $M_j$
are distinct Riemannian manifolds is much more difficult. 

Under Gromov-Hausdorff convergence 
\be
(X_j, d_j) \GHto (X_\infty, d_\infty),
\ee
one can apply Gromov's Embedding Theorem (\ref{Gromov-Z}) to define
what it means to say $x_j \to x_\infty$ for $x_j \in X_j$ 
as follows: there exists
a compact metric space, $Z$, and distance preserving maps, $\varphi_j: X_j \to Z$, such that
\be \label{Gromov-pts}
d_H^Z(\varphi_j(X_j), \varphi_\infty(X_\infty)) \to 0
\textrm{ and } d_Z(\varphi_j(x_j), \varphi_\infty(x_\infty))\to 0. 
\ee
Note that $x_\infty$ is not unique: if $F: X_\infty\to X_\infty$ is an
isometry then $F(x_\infty)$ is also a limit of $x_j$.  This is a natural consequence
of the fact that the GH distance is between isometry classes of metric
spaces.  Gromov shows that 
\be 
\forall \, x_\infty\in X_\infty \,\,\exists \,
x_j \in X_j \textrm{ such that } x_j \to x_\infty.
\ee
In fact, there exist functions, $H_j: X_\infty \,\to\, X_j$, and distance
preserving maps, $\varphi_j: X_j \to Z$, satisfying (\ref{Gromov-pts})
such that if $x_j=H_j(x_\infty)$ then $x_j \to x_\infty$ 
and
\be \label{dist-cont}
d_j(H_j(x), H_j(y)) \to d_\infty(x,y).
\ee
Furthermore, for any $r>0$, we have converging closed balls
\be  
x_j \to x_\infty \implies \left(\bar{B}(x_j, r), d_j\right) 
\GHto \left(\bar{B}(x_\infty, r), d_\infty\right).
\ee
By the compactness
of $Z$, there is a Bolzano-Weierstass Theorem: 
\be \label{GH-BW}
x_j \in X_j \implies \exists j_k \textrm{ s.t. } x_{j_k}\to x_\infty \in X_\infty.
\ee
Combining (\ref{dist-cont}) with (\ref{GH-BW}) we have
diameter continuity:
\be
\lim_{j\to\infty} \diam(X_j)=\diam(X_\infty).
\ee
See for example \cite{BBI} and \cite{Sormani-AA} for more details.

Now suppose we have intrinsic flat convergence,
\be\label{MjFtoM}
M_j = (X_j, d_j, T_j) \Fto M_\infty=(X_\infty, d_\infty, T_\infty).
\ee
In \cite{Sormani-AA}, the author applied the Intrinsic Flat Embedding Theorem
as in (\ref{converge}) to say a sequence $x_j \in X_j$ is {\em Cauchy}
if there exists
a complete metric space, $Z$, a point, $z_\infty \in Z$, and distance preserving maps, $\varphi_j: X_j \to Z$, such that
\be \label{point-conv-1}
d_F^Z\left(\varphi_{j\#} T_j ,\varphi_{\infty\#}T_\infty\right)\to 0
\ee
and $\varphi_j(x_j) \to z_\infty$.
One says the
sequence has {\em no limit} in $\bar{X}_\infty$ if 
\be\label{no limit}
z_\infty \notin \varphi_\infty(\bar{X}_\infty).
\ee
One says the points {\em converge} $x_j \to x_\infty$ in $X_\infty$ (or respectively
in $\bar{X}_\infty$) 
if there exists $x_\infty$ in $X_\infty$ (or respectively
in $\bar{X}_\infty$)
such that $z_\infty=\varphi_\infty(x_\infty)$:
\be\label{not-disappearing}
z_\infty \in \varphi_\infty(X_\infty)= \set( \varphi_{\infty\#}T_\infty)
\quad \left( \textrm{or respectively }
z_\infty \in \varphi_\infty(\bar{X}_\infty) \right).
\ee
Note that $x_\infty$ is not unique: if $F: X_\infty\to X_\infty$ is a current preserving
isometry then $F(x_\infty)$ is also a limit of $x_j$.  This is a natural consequence
of the fact that the $\mathcal{F}$ distance is between current preserving
isometry classes of integral current spaces. Below we will provide additional
conditions (\ref{balls-converge}) and (\ref{to-0}) which demonstrate that
Cauchy
sequences of points which disappear with respect to one 
sequence of $\varphi_j$ satisfying (\ref{point-conv-1})
cannot be Cauchy sequences of
points that converge with respect to another sequence
of $\varphi_j$ satisfying (\ref{point-conv-1})
and visa versa.

In \cite{Sormani-AA} the author proves that if (\ref{MjFtoM}) then
\be\label{to-a-point}
\forall \, x_\infty\in \bar{X}_\infty \,\,\exists \,
x_j \in X_j \textrm{ such that } x_j \to x_\infty.
\ee
In fact there exist {\em convergence functions}, $H_j: X_\infty \,\to\, X_j$, and distance preserving maps, $\varphi_j: X_j \to Z$, satisfying (\ref{point-conv-1})
such that $H_j(x_\infty) \to x_j$ as in (\ref{Gromov-pts})
with
\be \label{Hj}
d_j(H_j(x), H_j(y)) \to d_\infty(x,y) \textrm{ and }
H_j(\partial M_\infty) \subset \partial M_j.
\ee
In fact the author proves in \cite{Sormani-AA}[Theorem 5.1] the following:
\be
\textrm{ if }M_j \Fto M_\infty \textrm{ where } M_\infty \textrm{ is nonzero and precompact,}
\ee   
\be
\textrm{then there exists } N_j \subset M_j
\textrm{ such that } N_j \GHto \bar{M}_\infty.  
\ee
Additional
properties of the $N_j$ are described there.  In general intrinsic
flat limits need
not be precompact and need not have finite diameter \cite{SorWen2}.

Since $Z$ is only complete, we do not have a simple Bolzano-Weierstrass
Theorem.  In fact points may disappear under flat convergence even in
a compact $Z$.   So we do not have continuity of diameter.   However
we have semicontinuity of diameter
\be\label{diam-semi}
\liminf_{j\to \infty} \diam(M_j) \ge \diam(M_\infty).
\ee
and depth
\be\label{depth-semi}
\liminf_{j\to \infty} \depth(M_j) \ge \depth(M_\infty)
\ee
where 
\be\label{defn-depth}
\depth(M)=\sup\left\{ d_X(x,y): \,\, x\in X, \, y\in \set(\partial T)\,\right\} \in [0,\infty].
\ee

Recall that for almost every $r>0$ one may view a ball in an integral current
spaces as an integral current space itself
$
S(x,r)=\left(\,\set\left(T_j\rstr B(x,r)\right), \, d,\, T\rstr B\left(x,r\right)\right).
$
and examine the convergence of balls when $M_j \Fto M_\infty$.
In \cite{Sormani-AA} the author proves that if $x_j \to x_\infty$ as in (\ref{to-a-point}), then there is a subsequence $j_k$
such that 
\be  \label{balls-converge}
S(x_{j_k}, r) \Fto S(x_\infty, r) \neq {\bf{0}} \,\,\,\textrm{ {\em for almost every} } r>0.
\ee
If a Cauchy sequence of points, $x_j\in M_j$, has no limit in $\bar{X}_\infty$
as in (\ref{no limit}) then 
\be  \label{to-0}
\exists \delta>0 \,\,\,s.t.\,\,\,S(x_{j_k}, r) \Fto 0 \,\,\,\textrm{ {\em for almost every }} r\in (0,\delta).
\ee
Since (\ref{balls-converge}) and (\ref{to-0}) are intrinsic notions which do not
depend on a choice of distance preserving maps, one concludes that a Cauchy
sequence of points which has no limit in $\bar{X}_\infty$ with respect to one 
sequence of $\varphi_j$ satisfying (\ref{point-conv-1})
cannot be a Cauchy sequence of
points that converges with respect to another sequence
of $\varphi_j$ satisfying (\ref{point-conv-1})
and visa versa. Far more subtle is determining
when exactly a sequence of points converges to a point in 
$\bar{X}_\infty\setminus X_\infty$.


Let us further consider converging sequences of points, $x_j \to x_\infty$. 
As a consequence of mass semicontinuity as in (\ref{semicont}), we have
for almost every $r>0$
\be\label{semi-ball}
\liminf_{k\to \infty} \mass(S(x_{j_k}, r)) \ge \mass( S(x_\infty, r)).
\ee
As a consequence 
of (\ref{bndry-conv}), we have convergence of spheres
for almost every $r>0$
\be  \label{sphere-conv}
\partial S(x_{j_k}, r) \Fto \partial S(x_\infty, r) 
\ee
Combining sphere convergence as in (\ref{sphere-conv}) with 
filling volume continuity as in (\ref{fillvol-cont}), 
Portegies and the author have
proven that for almost every $r>0$
\be  \label{sphere-fillvol-cont}
\fillvol(\partial S(x_{j_k}, r)) \Fto \fillvol(\partial S(x_\infty, r))
\ee
as well as continuity of another notion called the sliced filling volume in \cite{Portegies-Sormani}.   In Theorem 4.27 of that paper, this continuity
is applied 
to determine when a sequence converges in $X_\infty=\set(T_\infty)$.  
The filling volume case of this theorem observes that
\be \label{conv-1}
\mass(S(x_\infty, r)) \ge \fillvol(\partial S(x_\infty, r)) = \lim_{k\to \infty}
\fillvol(\partial S(x_{j_k}, r)),
\ee
which implies that if there is a uniform lower bound $C>0$ such that
\be \label{conv-2}
\fillvol(\partial S(x_{j}, r))\ge Cr^m \qquad \forall j\in {\mathbb{N}}
\ee
then by
(\ref{def-set-current}), 
$x_\infty \in \set(T_\infty)=X_\infty $. 
An idea similar to this one was applied earlier in an extrinsic way by the 
author and Wenger in \cite{SorWen1} to show that when a sequence
of manifolds has nonnegative Ricci curvature or has a linear
contractibility function, then the GH and $\mathcal{F}$ limits
agree.  This method has been applied to prove intrinsic flat
and GH limits agree under a variety of different conditions by
the author with Portegies, by
Munn, by Perales, by Perales-Li,  and by Matveev-Portegies  
\cite{Portegies-Sormani} \cite{Munn-F=GH}   
 \cite{Perales-Conv} \cite{Li-Perales} \cite{Matveev-Portegies}.

With only lower scalar curvature bounds on the sequence one does not
in general have GH and intrinsic flat limits that agree (see examples
below) FILLIN and in fact one may not have any GH limit even for a
subsequence.   Nevertheless this method might be applied to determine
which points in the sequence of manifolds are disappearing and which
remain.   

In fact the author and Portegies prove a Bolzano-Weierstrass
Theorem for sequences of points with bounds on their filling or sliced filling volume
\cite{Portegies-Sormani}[Theorem 4.30].  One consequence of this
theorem is that if
a sequence satisfies (\ref{conv-1}) then a subsequence is Cauchy 
(using an argument involving the fact that the limit space has finite mass) 
and then by the above method, the sequence converges to a point in $X_\infty$.
There is also a Bolzano-Weierstrass Theorem for sequences of points 
which does not require a lower bound on the filling
volumes of spheres, but instead requires
\be
\exists r_0>0 \,\,\,s.t.\,\,\,
\liminf_{j\to\infty} d_{\mathcal{F}}(S(x_j,r), {\bf{0}}) \ge h(r)>0
\,\,\,\textrm{\em for a.e. } r\in (0, r_0]
\ee
to obtain a subsequence which is Cauchy and converges 
in $\bar{X}_\infty$ \cite{Sormani-AA}[Theorem 7.1].   Note that the 
theorems and proofs in \cite{Sormani-AA} are very easy to read.
The more technically difficult results in \cite{Portegies-Sormani} involving
filling volumes and sliced filling volumes are better in that they
provide more precise controls under weaker hypothesis.

\subsection{Limits of Functions}

Recall that when using the compactness and weak rigidity method to prove an almost rigidity theorem as in (\ref{alm-rig-Mj}), one proves that
$M_j \to M_0$, where $M_0$ is a specific given rigid space by first using a compactness theorem to show a subsequence, $M_{j_k} \to M_\infty$ and then proving $M_\infty=M_0$.   One way to prove that $M_\infty=M_0$ is to
construct an isometry between these spaces as a limit of functions from
$M_{j_k}$ to $M_0$. Theorems which produce limits of subsequences of 
functions are called Arzela-Ascoli Theorems. 

First let us describe what we mean by a limit function.  
Suppose $F_j : X_j \to Y_j$ are functions, and $X_j \to X_\infty$ and 
$Y_j \to Y_\infty$ in the GH or $\mathcal{F}$ sense then we say
$F_\infty: X_\infty\to Y_\infty$ is their limit, denoted
$F_j \to F_\infty$, if 
\be
F_j(x_j) \to F_\infty(x_\infty) \textrm{ whenever } x_j \to x_\infty
\ee
More precisely, $F_j \to F_\infty$ if there exists
convergence functions, $H_j: X_\infty \to X_j$ and $H'_j: Y_\infty\to Y_j$, 
as in (\ref{Hj}) such that
\be
F_j\circ H_j(x) = H_j' \circ F_\infty(x) \qquad \forall x \in X_\infty.  
\ee

The Gromov-Hausdorff Arzela-Ascoli Theorem states that if
one has compact metric spaces, $X_j \GHto X_\infty$
and $Y_j \GHto Y_\infty$, and if $F_j: X_j \to Y_j$ are equicontinuous
\be\label{equicont}
\forall \epsilon>0 \,\,\exists \delta_\epsilon>0\textrm{ such that }
d_{X_j}(x,x')< \delta_\epsilon \,\implies \,
d_{Y_j}(F_j(x), F_j(x'))\le \epsilon.
\ee
then $F_j \to F_\infty$ where $F_\infty: X_\infty\to Y_\infty$
satisfies (\ref{equicont}).  This theorem is a direct consequence of
Gromov's Embedding Theorem combined with the standard proof
of the classical Arzela-Ascoli Theorem (cf. \cite{Sormani-AA}).  Furthermore, if $F_j$ are surjective
then the limit $F_\infty$ is surjective.  If the $F_j$ are 
isometries on balls of radius $r_0>0$,
then the limit is an isometry on balls of radius $r_0>0$.  This was
applied by the author in joint work with Wei to prove the GH
limits of manifolds with $\Ricci\ge 0$ have
universal covering spaces \cite{SorWei1} and that the
covering spectrum is continuous with respect to Gromov-Hausdorff
convergence in \cite{SorWei3}.   

It is not absolutely necessary that the sequence of functions be
equicontinuous.  Gromov proves the same result for $\epsilon_j$
almost isometries $F_j: X_j \to Y_j$,
\be
|d_{Y_j}(F_j(p),F_j(q)) -d_{X_j}(p,q)|<\varepsilon_j \textrm{ and } 
Y_j \subset T_{\varepsilon_j}(X_j)
\ee
with $\epsilon_j \to 0$, producing a limit $F_\infty: X_\infty\to Y_\infty$
which is an isometry \cite{Gromov-metric}.  See also the Burago-Burago-Ivanov text 
\cite{BBI}.   Almost isometries are applied to prove
the GH almost rigidity theorems of Colding in \cite{Colding-volume} and 
of Cheeger-Colding in \cite{ChCo-almost-rigidity}
by constructing almost isometries from the $M_j$ to the rigid space,
$M_0$.  In \cite{Sor-cosmos}, a GH Arzela-Ascoli Theorem which only requires
almost equicontinuity is proven by the author in order to prove another
GH almost rigidity theorem. 



One cannot hope for a $\mathcal{F}$ Arzela-Ascoli Theorem which is as powerful as the GH Arzela-Ascoli Theorem.  In Example~\ref{ex-not-geod}
one has no limit for the geodesics $\gamma_j:[0,1] \to M_j$ running through
the increasingly thin tunnels.  Nevertheless two useful Arzela-Ascoli Theorems 
with additional hypotheses were proven  in \cite{Sormani-AA}.

Suppose $F_j: M_j \to M_j'$ are (surjective) current preserving isometries on balls
of radius $r_0>0$ and 
$
M_j \Fto M_\infty\neq {\bf{0}}$ and $M'_j \Fto M'_\infty\neq {\bf{0}}.
$
Then a subsequence of $F_j$ converges to $F_\infty: M_\infty\to M'_\infty$ which is also
a (surjective) current preserving isometry on balls
of radius $r_0>0$ \cite{Sormani-AA}.  This theorem has been applied
in joint work of the author with Sinaei to study the intrinsic flat convergence
of covering spaces and the covering spectrum \cite{Sinaei-Sormani-1}.

Suppose $F_j: M_j \to Y_j$ are equicontinuous maps as in (\ref{equicont})
where $M_j$ are integral current spaces and $Y_j$ are compact
metric spaces such that
$
M_j \Fto M_\infty\neq {\bf{0}}$  and $Y_j \GHto Y_\infty$.
Then $F_j$ converge to $F_\infty: M_\infty\to Y_\infty$ 
which also satisfies (\ref{equicont}).  Furthermore, if $F_j$ are surjective
then the limit $F_\infty$ is surjective \cite{Sormani-AA}.   
Keep in mind that
this includes equicontinuous functions, $F_j : M_j \to [a,b]$, and embeddings
in compact regions in Euclidean space, $F_j: M_j \to \mathbb{E}^N$.
This theorem has been applied jointly with Huang and Lee to prove
Almost Rigidity of the Positive Mass Theorem in the graph setting \cite{HLS}.

Conjectured related Arzela-Ascoli theorems
are suggested in \cite{Sormani-AA} and the 
proofs there are not difficult to read.  These theorems are particularly
useful when proving almost rigidity theorems to try to construct
isometries from the limit $M_\infty$ of a subsequence $M_j$ to the
desired rigid space, $M_0$.  

\begin{problem}
Suppose $F_j: M_j \to Y$ has $\Lip(F_j) \le 1$ with $Y$ compact
(including Riemannian embeddings
and graphs), and $M_j$ satisfy the hypothesis of
Wenger's Compactness Theorem as in (\ref{Wenger-comp}).
Then we know a subsequence 
$M_j \Fto M_\infty=(X_\infty, d_\infty, T_\infty)$.
We can see $F_{j\#}[M_j]$ satisfies the hypothesis of the Ambrosio-Kirchheim
Compactness as in (\ref{AK-compactness}),   
so we know a subsequence
$F_{j\#}[M_j]$ converges in the flat sense to some $S_\infty\in \intcurr_m(Y)$.
If $M_\infty\neq {\bf{0}}$
we know know by the Arzela-Ascoli Theorem above that a subsequence,
also denoted $F_j$, converges to $F_\infty: X_\infty \to Y$.
We conjecture that
\be \label{images-flat-converge}
M_\infty=(X_\infty, d_\infty, T_\infty) \textrm { with } F_{\infty\#} T_\infty=S_\infty.
\ee
\end{problem}

This conjecture captures the key steps used in joint work of the author
with Huang and Lee to prove the Almost Rigidity of
the Positive Mass Theorem in the graph setting in \cite{HLS}.
In that proof we have $M^m_j$ which satisfy our almost rigidity
conditions, and show they satisfy the hypothesis of Wenger's 
Compactness Theorem and then study their images as graphs in
$Y=\mathbb{E}^{m+1}$ and obtain (\ref{images-flat-converge}).
Some of the
techniques there may be useful towards proving this conjecture in
general.
  
\subsection{Intrinsic Flat with Volume Convergence}

Recall that intrinsic flat convergence does not imply
volume convergence. 
One does have semicontinuity, $M_j \Fto M_\infty$ implies 
$\liminf_{j\to \infty} \mass(M_j) \ge \mass(M_\infty)$, but even
when $M_j$ and $M_\infty$ are Riemannian manifolds the
volumes need not converge (as seen in the examples with
cancellation).  Note that all the examples presented
above with wells, 
bubbling and sewing do have volume convergence.
We introduce the following:

\begin{defn}\label{def-volF}
The intrinsic flat volume distance between two integral
current spaces, $M_j=(X_j, d_j, T_j)$ is
\be
d_{V\mathcal{F}}(M_1, M_2) = d_{\mathcal{F}}(M_1, M_2)
+|\mass(M_1) - \mass(M_2)|
\ee
\end{defn}

So $M_j \VolFto M_\infty$ iff $\mass(M_j) \to \mass(M_\infty)$
and $M_j \Fto M_\infty$.  Note that points may still disappear
as in the Ilmanen Example.   However
there is no cancellation.  

Portegies studied the properties of $M_j=(X_j, d_j, T_j)$ such that
$M_j \VolFto M_\infty$.
Applying his Lemma 2.7 in \cite{Portegies-F-evalue}, we see that the metric 
measure spaces, 
$(X_j, d_j, ||T_j||)$ converge in 
a measured sense.  That is, there are distance preserving maps, $\varphi_j: X_j \to Z$
such that 
\be\label{vol-to-measure}
\varphi_{j\*} ||T_j|| \to \varphi_{j\*} ||T_\infty|| \textrm{ weakly as measures in } Z.
\ee
Here
$Z$ is only complete and one need not have GH convergence (as in Ilmanen's
Example).   Portegies then applied this to prove that the Laplace eigenvalues of
converging sequences of manifolds are upper semicontinuous,
\be\label{evalues}
\limsup_{j\to \infty} \lambda_j(M_j) \le \lambda_j(M_\infty).
\ee
He presents examples showing this is false without 
$\mass(M_j) \to \mass(M_\infty)$.

Note that the Hawking Mass Compactness Conjecture and
Scalar Compactness Conjecture both propose that a subsequence 
$M_{j_k} \VolFto M_\infty$ as is shown in the compactness
theorem proven jointly with LeFloch \cite{LeFloch-Sormani-1}.  
Also note that in the work of the author with Lee, Huang, and Stavrov
proving various special cases towards almost rigidity of the Positive
Mass Theorem, it is proven that $M_j \VolFto M_0$ \cite{LeeSormani1}\cite{HLS}\cite{Sormani-Stavrov-1}.   Matveev and Portegies 
prove $M_{j_k} \VolFto M_\infty$ when Ricci curvature is uniformly
bounded below \cite{Matveev-Portegies}[Theorem 4.1].

Keep in mind that $M_{j_k} \VolFto M_\infty$ alone will not suffice
to achieve generalized $\Scal>0$ properties on $M_\infty$
as seen in the problematic examples involving sewing of a single 
curve to a point will have volume convergence and yet the
limit in (\ref{Scalar-volume}) will fail to be nonnegative for a ball
about that limit point.   Note however that if one assumes
(\ref{Scalar-volume-0}) holds at every point $p$ in every $M_j$ 
with a uniform lower bound on $r_p\ge r_0>0$, then (\ref{Scalar-volume-0})
can be shown to hold with inequalities on $M_\infty$.

\section{Results and Conjectures about Limits of Manifolds with Nonnegative Scalar Curvature}\label{sect-conj}

In this final section of the paper we state our almost rigidity conjectures
precisely and survey known results towards proving these conjectures.  We
suggest special cases which might be proven more easily.  Note that
completely proving any almost rigidity theorem is significantly more difficult that
proving the corresponding rigidity theorem.  One must either reprove
the rigidity theorem in a quantitative way obtaining (\ref{alm-rig-M}) using the
{\em explicit control method} as described in the introduction.  Or 
one must prove the Rigidity Theorem on a generalized class of spaces 
obtaining (\ref{alm-rig-Mj}) using the {\em compactness and weak rigidity method}.
Recall that we have already proposed two compactness
conjectures in Section~\ref{subsect-compactness}.    In the final two subsections
of this paper we discuss generalized notions of nonnegative
scalar curvature as proposed by Gromov and regularity theory.

\subsection{Almost Rigidity of the Positive Mass Theorem
and the Bartnik Conjecture}\label{sect-PMT}

Consider the class, $\mathcal{M}$, of asymptotically flat
three dimensional Riemannian manifolds with nonnegative
scalar curvature and no interior closed minimal surfaces
and either no boundary or the boundary is an outermost
minimizing surface.   This is the physically natural class of
spaces used to prove the Penrose Inequality as discussed
in the introduction.  By the Positive Mass Theorem, if $M\in \mathcal{M}$
has $m_{ADM}(M)=0$ then $M$ is isometric to Euclidean space.
 In \cite{LeeSormani1}, Lee and the
author proposed that almost rigidity for the Positive Mass Theorem
should be provable using intrinsic flat convergence and demonstrated
that it is false for GH convergence. 

\begin{conj}\label{conj-PMT}   
Fix $D>0$, $r_0>0$.
Let $M_j^3 \in \mathcal{M}$ and let $\Sigma_j\subset M_j^3$ be
special surfaces with $\area(\Sigma_j)=4\pi r_0^2$
and $\Sigma_\infty=\partial B(0, r_0)\subset {\mathbb{E}}^3$.  
We conjecture that
\be\label{LSconj}
m_{ADM}(M_j) \to 0 \implies d_{Vol\mathcal{F}}\left(\,T_D(\Sigma_j), T_D(\Sigma_\infty)\,\right)\to 0.
\ee
 where  
$T_D(\Sigma)$ is the 
tubular neighborhood of radius $D$ around $\Sigma$ or alternatively 
\be\label{HLSconj}
m_{ADM}(M_j) \to 0 \implies d_{Vol\mathcal{F}}\left(\,\Omega_j,\,\Omega_\infty\right)\to 0.
\ee
where $\Omega_j$ is the interior of $\Sigma_j$
with $\depth(\Omega_j)\le D$.  
 \end{conj}

 Note that Lee and
 the author were being deliberately vague as to what a special
 surface, $\Sigma_j$, should be in this conjecture.  A number of more precisely 
stated special cases of this 
conjecture were provided in the final section of \cite{LeeSormani1}  
along with brief ideas as to how
one might approach the proof of the conjecture in those cases.
Most of these special cases are still open.

 Lee and the author proved that (\ref{LSconj}) holds
 when $M_j^3$ have metric tensors
 of the form $g_j=dr^2 + f_j(r)^2 g_{{\mathbb{S}}^2}$
and $\Sigma_j=r^{-1}(r_0)$.   The proof uses 
Geroch monotonicity and Theorem~\ref{embed-const}
to obtain explicit controls on $T_D(\Sigma_j)$\cite{LeeSormani1}.
LeFloch and the author examined these explicit controls in more
detail in \cite{LeFloch-Sormani-1}.  The author and Stavrov
proved this conjecture when $M_j^3$ are Brill-Lindquist geometrostatic 
manifolds and $\Sigma_j$ are large spheres in \cite{Sormani-Stavrov-1}.  
That proof is also completed using explicit controls: 
bounds on the metric tensor are found on carefully selected regions 
within the manifolds followed by an application of Theorem~\ref{thm-subdiffeo}.

Huang, Lee and the author proved 
 (\ref{HLSconj}) when $\Omega_j^3\subset M_j$ are
graph manifolds, which have Riemannian embeddings, $\Psi_j: M^3_j \to 
{\mathbb{E}}^4$ as graphs.  We assumed controls on the $\Sigma_j$ and other technical properties on the 
$M_j^3$ \cite{HLS}.  Using the properties of
graph manifolds with $\Scal \ge 0$ and $m_{ADM}(M_j) \to 0$
we first proved that  $\vol(\Omega_j) \to \vol(B(0,r_0))$.
We applied Wenger's Compactness 
Theorem and an Arzela Ascoli Theorem
to prove a subsequence of the $\Omega_j \Fto \Omega_\infty$
and $\Psi_j$ converge to $\Psi_\infty: \Omega_\infty\to \mathbb{E}^4$
with $\Lip(\Psi_\infty)\le 1$.
By the lower semicontinuity of mass we showed
\be\label{massPsi}
\mass(\Psi_{\infty\#} [\Omega_\infty])
\le \mass( [\Omega_\infty])\le \vol(B(0,r_0)).
\ee
We used the controls on $\Sigma_j$ to prove they Lipschitz converge to
$\partial \Omega_\infty$ and that $\Psi_\infty: \partial \Omega_\infty \to 
\partial B(0, r_0)\times\{0\}$ is biLipschitz.   Thus
 $
 \partial \Psi_{\infty\#} [\Omega_\infty]=[\partial B(0,r_0) \times  \{0\}].
 $
Combining this with (\ref{massPsi}) implies $Psi_{\infty\#} [\Omega_\infty]$
solves the Plateau Problem.  So
\be
\Psi_{\infty\#} [\Omega_\infty]=\lbrack B(0, r_0) \times  \{0\}\rbrack
\textrm{ and }
\mass(\Psi_{\infty\#} [\Omega_\infty])=\vol(B(0,r_0)) (4/3).
\ee
This implies equality in (\ref{massPsi}) so 
$\Psi_\infty$ must be an isometry: $\Omega_\infty=B(0, r_0)$ \cite{HLS}.

\begin{rmrk}
One might consider applying the Huisken isoperimetric
mass as in (\ref{mISO}) to prove Conjecture~\ref{conj-PMT} with
$\Sigma_j$ chosen to be uniformly asymptotically spherical in $M_j$ with $m_{ADM}(M_j) \to 0$ so that 
\be \label{mISOto0}
m_{ISO}(\Omega_j)=\frac{2}{\area(\partial \Omega_j)}
\left(\vol(\Omega_j)-\frac{\area(\partial \Omega_j)^{3/2}}{6 \sqrt{\pi}}\right) \to 0.
\ee
\end{rmrk}

Observe that (\ref{mISOto0}) immediately implies
\be
\vol(\Omega_j)\to{A_0^{3/2}}/{(6 \sqrt{\pi})}= (4/3) \pi r_0^3=\vol(B(0,r_0))
\ee
So by Wenger's Compactness Theorem a subsequence converges
$\Omega_{j}\Fto \Omega_\infty$.  By lower semicontinuity
of mass we have 
$
\mass(\Omega_{\infty}) \le \vol(B(0,r_0)).
$

Suppose we also assume that there exist Riemannian isometric embeddings 
$\Psi_j: \Omega_j \to {\mathbb{E}}^N$ with the property that
when $\Psi_j$ restricted to $\partial \Omega_j$ is uniformly biLipschitz
to $\partial B(0,r_0)\times \{0,...,0\}$.  Then exactly as in the above description of the proof in \cite{HLS}
we have $\Omega_\infty$ isometric to $B(0,r_0)$.  Here the only place we
used $\Scal \ge 0$ was when we replaced $m_{ADM}(M_j) \to 0$
by (\ref{mISOto0}) considering that Miao's proof that $m_{ISO}(\Omega)$
is close to $m_{ADM}(M)$ for large round $\partial \Omega$ 
in \cite{Fan-Shi-Tam} uses $\Scal \ge 0$.  Note that $m_{ISO}(\Omega)$
for $M^3 \in \mathcal{M}$ does not imply $\Omega=B(0,r_0)$; one needs
to impose some asymptotic roundness on the $\partial \Omega$ 
as well as $\Scal \ge 0$ even to obtain rigidity. 

 \begin{problem} \label{BY-Compactness}
 Shi and Tam proved in \cite{Shi-Tam-2002} that the Brown-York Mass,
 $m_{BY}(\partial \Omega)$,
  is nonnegative if $\Scal \ge 0$ on $\Omega$
 and $\partial \Omega $ has positive Gauss curvature, and
  \be\label{BY-rigidity}
 m_{BY}(\partial \Omega)=0 \implies 
 \Omega \subset \mathbb{E}^3.
 \ee
This mass (which agrees with the Liu-Yau mass in this setting)
is defined using a Riemannian isometric embedding 
$\Psi: \partial \Omega \to \mathbb{E}^3$, the mean curvature, $H$,
of $\partial \Omega \subset M^3$ and the mean curvature, $H_0$, of
$\Psi(\partial \Omega) \subset {\mathbb{E}}^3$ as follows:
\be
m_{BY}(\partial \Omega)= \frac{1}{8\pi}\int_{\partial \Omega} H_0 - H \, d\sigma
\ee
The Arzela-Ascoli Theorems proven above might
be helpful towards proving Conjecture~\ref{Hawking-Compactness} for
the Brown-York Mass including
semicontinuity of the Brown-York mass under $\mathcal{F}$ convergence
and almost rigidity of the Shi-Tam Rigidity Theorem.   This might also
be applied to prove the almost rigidity of the Positive Mass Theorem
or the Bartnik Conjecture.  One of the biggest difficulties here is that 
mean curvature must be defined in a generalized way and controlled 
under intrinsic flat convergence.
\end{problem}

\begin{problem}\label{problem-HI}
Huisken and Ilmanen defined a weak notion of mean curvature
in their proof of the Penrose Inequality (which can also be applied
to prove the Positive Mass Theorem) \cite{Huisken-Ilmanen}.
Perhaps one might try to use their method to prove Almost Rigidity of
the Positive Mass Theorem.  One might consider a limit space
$M_\infty$ and attempt to define
Huisken-Ilmanen's weak inverse mean curvature flow and 
prove Geroch monotonicity on $M_\infty$.
What regularity is needed on $M_\infty$?  What notion of nonnegative
scalar curvature is required?
\end{problem}

Recall that in \cite{Geroch-monotonicity}, Geroch  proved that if $N_t: \mathbb{S}^2 \to M^3$
evolves by smooth inverse mean curvature flow,
\begin{eqnarray}\label{IMCF}
\frac{d}{dt} x = \frac{\nu}{H}
&\textrm{ where }&\,\, x=N_t(p), \,\, \nu\textrm{ is the normal to } N_t
\textrm{ at } x,\\
&& \textrm{ and } H \textrm{ is the mean curvature of }
N_t \textrm{ at } x,
\end{eqnarray}
and $M^3$ has $\Scal \ge 0$ then
then the Hawking mass, $m_H(N_t)$, is nondecreasing. 
\be \label{G-mono}
t_2 > t_1 \implies m_H(N_{t_2}) \ge m_H(N_{t_1}).
\ee
In \cite{Huisken-Ilmanen},
Huisken-Ilmanen introduced weak inverse mean curvature flow,
\be\label{wIMCF}
N_t=\partial \{x:\,\, u(x)<t\} \,\,\,\textrm{ where } \,\,\,
\textrm{div}_M\left(\frac{\nabla u}{|\nabla u|} \right)=|\nabla u|,
\ee
proving it also satisfies Geroch monotonicity 
as in (\ref{G-mono}), 
 and 
with the right boundary conditions 
$\lim_{t\to\infty} m_H(N_t)= m_{ADM}(M)$.   They defined a weak
mean curvature for the level sets of $u$ to be $H=\nabla u$ almost
everywhere.   One may naturally ask what regularity was needed on
the limit space to prove these results.  

Alternatively one might apply Huisken-Ilmanen's method on the sequence
$M_j$ rather than on $M_\infty$.  One might consider sequences 
of $u=u_j$ satisfying (\ref{wIMCF}) on $M_j$
and consider limits $u_j \to u_\infty$, if it is difficult to define 
(\ref{wIMCF}) on $M_\infty$ itself.   On each $M_j$ one can define (\ref{wIMCF})
and then one has $m_H(N_t)$ nearly constant.  This was a key step in the proof of
the Almost Rigidity of the Positive Mass Theorem in the rotationally
symmetric case \cite{LeeSormani1}.  It would be interesting to investigate this
even in
the setting with smooth inverse mean curvature flow.  In the setting where
one only has weak inverse mean curvature flow, then the $N_t$ may
skip over entire regions.  In private communication with the author,
Huisken has suggested that it
might be possible to bound the volume of the skipped regions using
his isoperimetric masses.   These same techniques might also be applied to
prove the Almost Rigidity of the Penrose Inequality.

\begin{problem}
Recall that the Penrose
Inequality for $M^3\in \mathcal{M}$: 
\be
m_{ADM}(M^3)\ge m_H(\partial M^3) = \sqrt{\tfrac{\area(\Sigma)}{16\pi}}
\ee
and Penrose Rigidity:
\be
m_{ADM}(M^3)= m_H(\partial M^3) \implies M^3 \textrm{ is isometric to }M_{Sch,m}
\ee
where $M_{Sch,m}$ is the Riemannian Schwarschild space with
mass $m=m_{ADM}(M^3)$:
\be
M_{Sch,m}=\left({\mathbb{R}}^3\setminus \{0\}, (1+ \tfrac{m}{2r})^2 \delta\right).
\ee
Almost rigidity for the Penrose Inequality was conjectured jointly with
Lee in \cite{LeeSormani2}  where it was proven
in the rotationally symmetric setting.   This has not yet been explored for
graph manifolds nor for Brill-Lindquist Geometrostatic manifolds.   These 
are perhaps easy enough to assign as a first project to a doctoral student
as the techniques in \cite{HLS} and \cite{Sormani-Stavrov-1} should directly
apply.  Proving it in general would involve all the same difficultes
as proving the almost rigidity for the Positive Mass Theorem and more.
\end{problem}

\subsection{The Bartnik Conjecture}

Bartnik's quasilocal mass of a region $\Omega_0 \subset M^3$
where $M_0^3\in \mathcal{M}$ and $\partial M_0^3 \subset \Omega_0$
was defined in \cite{Bartnik-1986} as an infimum of
the ADM masses of extensions, $M$, of $\Omega_0$:  
\be
m_B(\Omega_0)=\inf\{m_{ADM}(M)|\, M \in\mathcal{PM}(\Omega_0) \}
\ee
where $M\in \mathcal{PM}\subset \mathcal{M}$ if it contains an isometric image
of $\Omega$: $\Omega \subset M$.   Bartnik conjectured that this
infimum is achieved by what he called the {\em minimal mass extension}
and that this minimal mass extension is scalar flat and static.   
Significant
research in this area has been completed by 
Corvino in \cite{Corvino-CMP}
in which the properties of a minimal mass extension are proven assuming it
exists.  Miao has searched for static extensions using perturbative methods
in \cite{Miao-CMP}.
To prove the minimal mass extension exists in general one may consider the
following approach:

\begin{conj}[Bartnik Conjecture]
There exists a sequence
\be\label{Bartnik-sequence}
M_j \in\mathcal{PM}(\Omega_0) \textrm{ such that }m_{ADM}(M_j) \to m_B(\Omega_0).
\ee
with a limit $M_j \to M_\infty$  such that
$
m_{ADM}(M_\infty)=m_B(\Omega_0)=\lim_{j\to \infty} m_{ADM}(M_j),
$
where $M_\infty$ is asymptotically flat, smooth, scalar flat and static.
\end{conj}

By the definition of Bartnik Mass we know there is a sequence
satisfying (\ref{Bartnik-sequence}).   We propose that one 
may be able to prove an intrinsic
flat compactness theorem for $\Omega_{R,j} \subset M_j$ 
where
$\Omega_{R,j}=T_R(\partial \Omega)$ or perhaps $\Omega_{R,j}\supset \Omega_0$ 
with well chosen $\partial \Omega_{R,j}=\Sigma_{R,j}\subset M_j$.   Note that by Wenger's
Compactness Theorem we only need 
\be
\vol(\Omega_{R,j})\le V_R \textrm{ and }
\vol(\partial \Omega_{R,j})\le A_R
\ee
to obtain a subsequence $\Omega_{R,j} \to \Omega_{R,\infty}$
since
$
\diam(\Omega_{R,j})\le 2R + \diam(\Omega_0).
$
Then we could glue together an $M_\infty$ from these limit regions
$\Omega_{R,\infty}$ and $\Omega_0\subset M_\infty$.   Our conjectured
Hawking mass compactness theorem would imply that $M_\infty$
has generalized $\Scal\ge 0$ and a well controlled Hawking mass.

However one must be warned that an $M_\infty$ obtained in this manner
need not be asymptotically flat.   In \cite{LeeSormani2},
Lee and the author proved that a sequence $M_j$ approaching equality in the
Penrose Inequality may develop a longer and longer neck.  The recent examples of Mantoulidis and Schoen which satisfy (\ref{Bartnik-sequence}) also develop
increasingly long necks so that $M_j$ smoothly converge on regions around $\Omega_0$ to an $M_\infty$ which is not asymptotically flat and has no ADM mass \cite{Mantoulidis-Schoen-16}.  For the regions near infinity these $M_j$ are 
simply Schwarzschild space.  Thus a new construction is needed which shortens
these necks so that the $M_j$ are in some sense uniformly asymptotically flat.
Then one could try to prove a subsequence converges in the intrinsic flat sense 
to an aymptotically flat limit.

Once one has proven $M_j \to M_\infty$ where $M_\infty$ is
asymptotically flat and has generalized $\Scal\ge 0$, then one must prove
the semicontinuity of the ADM masses.
In \cite{Jauregui} 
Jauregui has proven this semicontinuity in a variety of settings including intrinsic flat
convergence in the rotationally symmetric case applying the Hawking mass
Compactness Theorem proven by the author with LeFloch in \cite{LeFloch-Sormani-1}.   
Jauregui-Lee have proven semicontinuity of the ADM mass
under $C_0$ convergence using the Huisken isoperimetric mass
in \cite{Jauregui-Lee:Huisken-isoper}.   These techniques might apply more generally.



\subsection{Almost Rigidity of the Scalar Torus Theorem}
\label{sect-torus}

The Scalar Torus Theorem states any $M^n$ diffeomorphic
to ${\mathbb{T}}^n$ with $\Scal \ge 0$ is isometric to a flat torus.  It
was proven using minimal surfaces for $n\le 7$ by Schoen and Yau in
\cite{Schoen-Yau-min-surf}.  It was proven by Gromov and
Lawson in all dimensions using spinors and the Lichnerowetz formula in  \cite{Gromov-Lawson-1980}.
Gromov proposed the following almost rigidity
conjecture vaguely in \cite{Gromov-Dirac}.    To avoid collapsing and expanding we have normalized the manifolds with the volume and diameter bounds.
To avoid bubbling, we have added the uniform lower bound on $\mina$
of (\ref{mina}).
This conjecture is false for GH convergence as seen in Example~\ref{well-tori}.

\begin{conj}
Let $M_j^3$ be diffeomorphic to $\mathbb{T}^3$
with  $\Scal(M_j^3) \ge -1/j$ and 
\be
\vol(M_j^3)=V_0 \quad \diam(M_j^3)\le D_0 \quad \mina(M_j^3)\ge A_0>0.   
\ee
Then a subsequence $M_{j_k} \Fto M_0$ where $M_0$ is a flat torus.
Possibly $M_{j_k} \VolFto M_0$.
\end{conj}

Note that by Wenger's Compactness Theorem, we know $M_{j_k} \Fto M_\infty$,
and by semicontinuity we have $\diam(M_\infty) \le D_0$ and 
$\mass(M_\infty)\le V_0$.   We do not know if $M_\infty$ is the ${\bf{0}}$
integral current space.  So even proving that much would be interesting.
Some points will disappear as seen in Example~\ref{well-tori}, so one must find
a sequence of points which doesn't disappear.

Even assuming that $M_\infty\neq {\bf{0}}$ one needs to prove that it has
some sort of generalized $\Scal \ge 0$ which is strong enough to prove torus rigidity.   Gromov suggests a few such notions which should be strong enough
in \cite{Gromov-Dirac} and prove $C^0$ limits of $M^3_j$ as above satisfy
these conditions. Examining his proofs and considering filling volumes 
and sliced filling as defined
in joint work of the author with Portegies \cite{Portegies-Sormani} might be helpful.
Note that Bamler has proven $\Scal \ge 0$ persists under $C^0$
convergence using Ricci flow in \cite{Bamler-16}.

Perhaps a more approachable problem would be to consider first $M^3_j$
which are graphs over the standard ${\mathbb{T}}^3$ and try to prove
$M^3_j \Fto {\mathbb{T}}^3$ using methods similar to those used in joint work
of the author with Huang and Lee \cite{HLS}.   Another possibility is to consider $M_j^3$ with metric tensors $g_j=dx^2+dy^2 +f_j(x,y) dz^2$ or $g_j=a_j(z)^2dx^2+b_j^2(z)dy^2 + dz^2$ and first try to prove
the warping functions $a_j, b_j$ and $f_j$ converge in the $H^1_{loc}$ sense 
to a metric with generalized $\Scal\ge 0$ as in joint work of the 
author with LeFloch \cite{LeFloch-Sormani-1}.  It is possible that spinors and
the Lichnerowetz formula might be formulated weakly for an $H^1_{loc}$
metric tensor.  More likely one can find a partial differential inequality 
on the warping functions which implies the rigidity theorem and can be shown 
to persist weakly under $H^1_{loc}$.  Finally one can use an Arzela Ascoli Theorem to relate the limit space obtained under $H^1_{loc}$ convergence and the
intrinsic flat limit.

\subsection{Almost Rigidity Theorems and Ricci Flow}

Bray, Brendle and Neves proved the
Cover Splitting
Rigidity Theorem as in (\ref{Cover-Splitting-Rigidity}) and
Bray, Brendle, Eichmair and Neves have proven the
$\mathbb{RP}^3$ Rigidity
Theorem as in (\ref{RP3-Rigidity}) using Ricci flow \cite{BBN-CAG-10}
\cite{BBEN-CPAM-10}.    Here we propose
the following:

\begin{conj} [Almost Rigidity of the ${\mathbb{RP}}^3$ Rigidity Theorem]
$\qquad\qquad\qquad\qquad$
Given a sequence of $M_j^3$ that are diffeomorphic to ${\mathbb{RP}}^3$ with
\be
\Scal_j \ge \tfrac{6j}{j-1} \quad \mina(M^3)= 2\pi \quad
\vol(M_j^3)\le V_0 \quad \diam(M_j^3)\le D_0
\ee
then $M_{j}^3 \Fto M_0=\mathbb{RP}^3$ or 
possibly even $M_{j}^3 \VolFto M_0=\mathbb{RP}^3$.
\end{conj}

One may construct examples with wells demonstrating
that $M_{j}^3 \GHto \mathbb{RP}^3$ may fail.

\begin{conj} [Almost Rigidity of Cover Splitting Theorem]
$\qquad\qquad\qquad\qquad$
Given a sequence of $M_j^3$ that contain noncontractible ${\mathbb{S}}^2$
and have
\be 
\Scal_j \ge \tfrac{2j}{j-1} \quad \mina(M^3)=4\pi \quad
\vol(M_j^3)\le V_0 \quad \diam(M_j^3)\le D_0
\ee
then a subsequence $M_{j_k}^3 \Fto M_0$
where 
$\tilde{M}_0$ is isometric to $\mathbb{S}^2\times \mathbb{R}$.
\end{conj}

\begin{problem}\label{well-cylinder}
Can one construct $M_j$ with $Scalar \ge 2j/(j-1)$ that are diffeomorphic
to ${\mathbb{S}}^2 \times {\mathbb{S}}^1$
and contain balls of radius $1/4$ that are isometric to 
balls in rescaled standard spheres?  If so then one can attach wells
and have no GH limit.
\end{problem}

One approach to proving these conjectures would be to use 
volume renormalized Ricci flow, $M_{j,t}$, of $M_j$ and show 
$M_{j,t}$ flows as $t\to \infty$ to an $M_{j,\infty}$ which is
isometric to $M_0$.  One may simplify things by assuming
a smooth Ricci flow exists for all time, or try to prove this,
or attempt to deal with Ricci flow through singularities.  If
\be\label{Ricci-flow}
d_{Vol\mathcal{F}}(M_{j,t}, M_{j,s}) < \epsilon\left(\tfrac{1}{t}-\tfrac{1}{s}\right)
\textrm{ where } \lim_{\delta\to 0} \epsilon(\delta) =0
\ee   
where $\epsilon$ is independant of $j$, then we have the conjecture.
Continuity of Ricci flow with respect to the intrinsic flat distance has
been studied by Lakzian in \cite{Lakzian-Ricci-flow} but only to analyze the
Ricci flow through a neck pinch singularity.  The author believes one might
be able to obtain (\ref{Ricci-flow}) constructing an explicit 
\be
Z=\{(x,r)|\, x\in M_{j,(1/r)},\, r\in [(1/t), (1/s)]\} 
\textrm{ with } g=dr^2 + f(r)^2 g_{j,(1/r)}.
\ee 
Again one might consider special cases where $M_j$ are known to be warped products.

\subsection{Gromov's Prisms and Generalized Nonnegative Scalar Curvature} \label{sect-generalized}

Gromov suggested in \cite{Gromov-Plateau} that if a sequence of $M_j \Fto M_\infty$ has $\Scal\ge 0$ then $M_\infty$ might have $\Scal\ge 0$ in some generalized sense.  In \cite{Gromov-Dirac}, Gromov proved the Gauss Bonnet
Prism Inequality for prisms in manifolds with $\Scal\ge 0$: 
\be\label{Prism}
\sum_{i=1}^3 \alpha_i \ge \pi  \textrm{ and } \sum_{i=1}^3 \alpha_i = \pi\implies
P \textrm{ is a prism in } \mathbb{E}^3
\ee
where $\alpha_i$ are bounds on the dihedral angles between the sides of the
prism which are minimal surfaces.   He suggests that the Prism Inequality
could be used to define generalized $\Scal\ge 0$ and then be applied to prove the torus
rigidity theorem on such limit spaces.  

Before we can define such a generalized notion of nonnegative scalar curvature
on integral current spaces, we would need to prove more regularity on the
$\mathcal{F}$ limit spaces.  Ordinarily there is no notion of angle on an integral
current space, nor is an integral current space a geodesic space.  In fact it might
even be the zero space. It should be noted that Ilmanen, the author and Wenger
conjectured in \cite{SorWen2}[Conjecture 4.18] that ``a converging sequence of three dimensional Riemannian manifolds with positive scalar curvature, a uniform lower bound on volume, and no interior closed minimal surfaces converges 
without cancellation to a nonzero integral current space.''  Combining Gromov's
ideas with this one we introduce the following conjecture:

\begin{conj}\label{fund-conj}
Suppose $M^3_j \Fto M_\infty$ or $M^3_j \VolFto M_\infty$ where
 $\Scal_j\ge 0$ and
\be
 \mina(M_j)\ge A_0>0 \quad \vol(M_j) \in [V_0,V_1]\subset (0,\infty) \quad \diam(M_j)\le D_0.
\ee
Then we have the following:

\vspace{.2cm}
\noindent
(a) $M_\infty$ is a nonzero integral current space.   In fact there is no
cancellation without collapse: If $p_j\in M_j$ has no limit in $\bar{M}_\infty$ then $\exists \delta>0$
s.t. $\vol(B(p_j, \delta)) \to 0$.

\vspace{.2cm}
\noindent
(b) $M_\infty$ is geodesic: If $p,q\in M_\infty$ there exists $p_j \to p$
and $q_j \to q$ with midpoints $x_j$ that converge to $x_\infty$ which
is a midpoint between $p$ and $q$ in $M_\infty$.

\vspace{.2cm}
\noindent
(c) There is a notion of angle between geodesics emanating from a point.

\vspace{.2cm}
\noindent
(d) There is a notion of dihedral angle
between two surfaces at $p\in \Sigma\cap \Sigma'\subset M_\infty$.

\vspace{.2cm}
\noindent
(e) Gromov's Gauss-Bonnet Prism Inequality as in (\ref{Prism}) holds on $M_\infty$.

\vspace{.2cm}
\noindent
(f) $\forall p\in M_\infty \,\, \exists r_p>0 \,\,s.t. \,\, \forall r<r_p\,\,\,\, 
V_p(r)=\mass(B(p,r)) / (4\pi r^3/3)\le 1.$

\end{conj}

Note that (a)-(d) are regularity properties for our limit spaces.
As seen in Examples above, they do not hold on limits of manifolds
with $\Scal_j\ge 0$ unless the $\mina(M_j)\ge A_0>0$.   The notion
of sliced filling volume developed in joint work with Portegies in \cite{Portegies-Sormani} might
be helpful towards proving (a) and (b).   It is quite possible that (c) is false
but that one can still prove (d) using the limit process to define the
dihedral angle.   In the special case where there exists Riemannian embeddings
$\Psi_j: M_j \to \mathbb{E}^N$, one can use an Arzela-Ascoli Theorem
to obtain $\Psi_\infty: M_j \to \mathbb{E}^N$ and then use $\mathbb{E}^N$
to define angles as needed in (c) or (d) and examine the semicontinuity of such angles.  
Note that any Riemannian manifold satisfies (a)-(d)
so they do not capture a generalized notion of nonnegative
scalar curvature.  

Now (e)-(f) are properties which capture $\Scal\ge 0$.  Property (e) proposed by Gromov needs a notion of angle so one needs to prove (d) first or use embeddings
into some large $\mathbb{E}^N$ to even define what this means.  Then one need only prove semicontinuity of the angles.  The power of property (e) is described
in \cite{Gromov-Dirac} including ideas towards the possibility that (e) implies (f).  Recall that (f) on a Riemannian manifold is equivalent to $\Scal \ge 0$ but was not powerful enough to prove any global properties of such spaces.  Nor is (f) continuous with respect
to $Vol\mathcal{F}$ convergence unless one assumes a uniform lower bound $r_p\ge r_0>0$ for all $p\in M_j$ for all $j \in {\mathbb{N}}$.    Nevertheless any
natural notion of $\Scal \ge 0$ on an integral current space ought to imply (f).
Note that Gromov also proves hyperbolic and spherical prism
inequalities on spaces with lower bounds on scalar curvature which are negative or positive respectively and that such inequalities might be used to generalize these
lower bounds on scalar curvature and one might prove they imply
the appropriate limit as in (\ref{Scalar-volume}) on $M_\infty$.
 
 \begin{rmrk}
In Remark~\ref{BY-Compactness}
we proposed that the Brown York mass might
be semicontinuous with respect to $\mathcal{F}$ convergence.
One might devise a way to define the Brown-York 
mass, $m_{BY}(\partial P)$, where $P$ is a prism.   These definitions
are likely only to involve integrals of dihedral angles.   It is quite possible that it would be easier to study the limits of Brown-York
masses of prisms than to even define dihedral angles on the limit spaces.  
If one can do this, then might try to replace $(e)$ with 

\vspace{.2cm}
\noindent
$(e_{BY})$  Shi-Tam Nonnegativity of the Brown-York Mass of Prisms: $m_{BY}(\partial P)\ge 0$.
\vspace{.2cm}

\noindent
Perhaps
the consequences Gromov devises using (e) might be concluded from 
$(e_{BY})$.
\end{rmrk}

\subsection{Almost Rigidity of the Positive Mass Theorem and Regularity
of $\mathcal{F}$ Limits}

It should be noted that in the famous work of Cheeger-Colding on the properties
of metric measure limits of Riemannian manifolds with lower bounds on 
Ricci curvature, a key step in proving regularity and the existence of
Euclidean tangent spaces at regular points, was the proof of their Almost Splitting
Theorem \cite{ChCo-PartI}.  In their work a point $p$ in a limit space $M_\infty$ 
has a Euclidean tangent space if the sequence of
rescaled balls
\be
(B(p,r_j), d_j/r_j)\textrm{ where } r_j \to 0
\ee
converges in the GH sense to a ball $B(0,1)$ in Euclidean space.
The Toponogov Splitting Theorem 
was similarly used to prove regularity results for
Alexandrov Spaces by Burago-Gromov-Perelman \cite{BGP}. 
However there is no splitting theorem for manifolds with 
nonnegative scalar curvature.

\begin{rmrk}\label{rmrk-regularity}
Let us consider $M_j \Fto M_\infty$ where the $M_j$ satisfy conditions
like those in the proposed Scalar Compactness Conjecture,
\be\label{scalar-mina-regularity}
\Scal_j \ge 0 \textrm{ and } \mina(M_j) \ge A_0>0
\ee
 and $p\in M_\infty$.
If one has proven a local version of the almost rigidity of the Positive Mass
Theorem then one can determine settings in which 
\be\label{regular-point}
(B(p,r_j), d_j/r_j, [B(p,r_j)]) \Fto (B(0,1), d_{{\mathbb{E}}^3}, [B(0,1)]).
\ee
This may possibly then be applied in the place of a splitting theorem
to prove some sort of regularity on the limit space, $M_\infty$.  
\end{rmrk}

Note that without the assumptions in (\ref{scalar-mina-regularity}), an intrinsic
flat limit may have no regular points in the sense described in
(\ref{regular-point}).   In the original paper defining intrinsic flat convergence
by the author with Wenger \cite{SorWen2}, an example of a sequence of
$M_j$ is given which converges in the intrinsic flat sense to taxicab space,
\be
M^2_{taxi}=([0,1]\times[0,1], d_{taxi}, [ \, [0,1]\times[0,1] \, ])
\ee where
\be
d_{taxi}\left((x_1,x_2), (y_1, y_2)\right)=|x_1-y_1|+|x_2-y_2|.
\ee
Such a space has no notion of angles and no prism properties.   It is possible
to imagine how sewing methods could be used to construct a sequence of
$M_j^3$ with $\Scal_j \ge -1/j$ that converges to $M^3_{taxi}$.   However, the $\mina(M_j) \ge A_0>0$ condition fails on such examples.  The taxi limit example
in \cite{SorWen2} satisfies $\mina(M_j) \ge A_0>0$ but constains points for
which the scalar curvature decreases to negative infinity.

Given $M_\infty$ as in Remark~\ref{rmrk-regularity}
one may even be able to prove that almost every point $p\in M_\infty$
satisfies (\ref{regular-point}).   Thus it may well be worthwhile to examine
to what extent one may use this kind of regularity to define dihedral angles,
mean curvatures, quasilocal masses and
a generalized notion of nonnegative scalar curvature on $M_\infty$.  

\vspace{.5cm}

\noindent
{\bf Acknowledgements:}
   The author would like to
thank the following mathematicians for discussions over the past few years related to this paper and for suggested references: Bartnik, Brendle, Carlotto, Cederbaum, Eichmair, Gigli, Huisken, Ilmanen, Khuri, Gromov, Huang, Jauregui, Lawson, Lee, LeFloch, Miao, Morgan, Sakovich, Schoen, Shi, Stavrov, Wang and Yau.  The author's research was funded by an individual research grant NSF-DMS-1309360 and a PSC-CUNY Research Grant.  This survey was written as notes for a course at the 2016 School in Geometric Analysis at the Lake Como School for Advanced Studies sponsored by IISS, ERC, GNAMPA, and DISAT, organized by Besson, Pigola, Impera, Rimoldi, Setti, Troyanov, Valtorta, and Veronelli.


\bibliographystyle{alpha}
\bibliography{2016-RR}

\end{document}